\documentclass[11pt]{article}\usepackage[dvips]{graphicx}

\def\mockpicture#1{%
\vbox to 130pt{\hsize 130pt\hrule
\noindent\vrule height130pt\hfill\vrule\hrule}}

\usepackage{amsmath}
\usepackage{amsfonts}

\newcommand{\A}{{\cal A}}
\newcommand{\B}{{\cal B}}
\newcommand{\C}{ {\mathbb{C}} }
\newcommand{\CC}{{\bf C}}
\newcommand{\Cl}{{\bf C}}
\newcommand{\E}{{\bf E}}
\newcommand{\EE}{{\cal E}}
\newcommand{\F}{{\cal F}}
\newcommand{\G}{{\cal G}}
\newcommand{\LL}{{\cal L}}
\newcommand{\PP}{{\cal P}}
\newcommand{\R}{ {\mathbb{R}} }
\newcommand{\GS}{{\cal S}}
\newcommand{\GSn}{{\cal S} ( \, [n] \, ) }
\newcommand{\GSpq}{{\cal S} ( \, [p+q] \, ) }
\newcommand{\U}{{\cal U}}
\newcommand{\Z}{ {\mathbb{Z}} }

\newcommand{\sanc}{{\cal S}_{\hbox{{\scriptsize\it ann-nc}}} (p,q) }
\newcommand{\stanc}{{\cal S}_{\hbox{{\scriptsize\it ann-nc}}} (2p,2q)}
\newcommand{\sancpair}{{\cal S}_{\hbox{{\scriptsize\it ann-nc-pair}}}(2p,2q) }
\newcommand{\smanc}{{\cal S}_{\hbox{{\scriptsize\it ann-nc}}} (p_1, \ldots , p_l) }
\newcommand{\sdnc}{{\cal S}_{\hbox{{\scriptsize\it disc-nc}}} (n) }

\newcommand{\card}{ \mbox{card} }
\newcommand{\perm}{perm}
\newcommand{\orbits}{ \mbox{orbits} }
\newcommand{\tr}{ \mbox{tr} }
\newcommand{\ncdisc}{ NC_{disc} }
\newcommand{\ncannpq}{ NC_{ann} (p,q) }

\makeatletter\def\footnotemark{\@ifnextchar[{\@xfootnotemark}
{\stepcounter{footnote}\begingroup\let\protect\noexpand
\xdef\@thefnmark{(\thefootnote)}\endgroup\@footnotemark}}
\makeatother

\font\titlefont=ptmri at 20pt

\def\capt#1{\vbox{\hsize150pt\raggedright\noindent{\small#1}}}
\def\doublecapt#1{\vbox{\hsize300pt\raggedright\noindent{\small#1}}}

\begin{document}

\title{\titlefont 
Annular non-crossing permutations and partitions, and 
second-order asymptotics for random matrices}

\author{{\sc James A. Mingo}\thanks{ 
Research supported by a grant from the  Natural Sciences and
Engineering Research Council, Canada.} \ \
{\sc and Alexandru Nica}$^{(\ast)}$\thanks{
Research supported by a Premier's Research Excellence Award,
Ontario, Canada.}}

\date{ }
\maketitle

\begin{abstract}
We study the set $\sanc$ of permutations of $\{ 1, \ldots, p+q \}$ which are 
non-crossing in an annulus with $p$ points marked on its external circle and
$q$ points marked on its internal circle. We define $\sanc$ algebraically, by 
identifying the crossing patterns which can occur in an annulus. We prove the 
annular counterpart for a ``geodesic condition'' shown by Biane to 
characterize non-crossing permutations in a disc. We examine the relation 
between $\sanc$ and the set $\ncannpq$ of annular non-crossing partitions of 
$\{ 1, \ldots , p+q \}$, and observe that (unlike in the disc case) the 
natural map from $\sanc$ onto $\ncannpq$ has a pathology which prevents 
it from being injective.

\vspace{6pt}

We point out that annular non-crossing permutations appear in the description 
of the second order asymptotics for the joint moments of certain families 
(Wishart and GUE) of random matrices. Some of the formulas extend to a 
multi-annular framework; as an application of that, we observe a phenomenon 
of asymptotic Gaussianity for traces of words made with independent 
Wishart matrices.
\end{abstract}

\section{Introduction} 

The set $NC(n)$ of non-crossing partitions of $\{ 1, \ldots , n \}$ has
\setcounter{page}{1}
been studied as an important example of a lattice, at least since the work
of Kreweras \cite{K}. It is customary (also since \cite{K}) to draw these 
partitions as pictures in a disc with $n$ points marked around its boundary.
In order to emphasize that, we will use in this paper the notation 
$\ncdisc (n)$ (instead of just $NC(n)$).

In work related to combinatorial aspects of non-commutative probability,
Biane \cite{B1}, \cite{B2} noticed that it is advantageous to embed 
$\ncdisc (n)$ into the group $\GS (n)$ of
permutations of $\{ 1, \ldots , n \}$; we thus arrive to talk about the
the set of disc non-crossing {\em permutations,} $\sdnc$. While $\sdnc$ 
is ``almost the same thing'' as $\ncdisc (n),$ it is nevertheless worth
being considered, because it has a non-trivial
equivalent characterization in terms of a ``geodesic condition'' in the
Cayley graph of $\GS (n)$ (cf. the review of this done in Section 2.10 
below).

\vspace{10pt}

In this paper we study the set $\sanc$ of permutations of 
$\{ 1, \ldots ,p+q \}$ which are non-crossing in an annulus 
with $p$ points marked on its external circle and $q$ points marked on
its internal circle.
 
\vspace{10pt}

The first thing to clarify is what is the formal algebraic definition 
of $\sanc$. This has to be an analogue for the fact that a permutation 
$\tau$ of $\{ 1, \ldots , n \}$ is in $\sdnc$ precisely when every cycle 
of $\tau$ is ``standard'' (in a natural sense), and when $\tau$ does not
display the crossing pattern ``$(a,c)(b,d)$'' with $1 \leq a<b<c<d \leq n$
(cf. the review done in Sections 2.4--2.8).
It turns out that in the annular framework there are three possible 
crossing patterns, rather than just one -- see the conditions (AC-1), 
(AC-2), (AC-3) in Definition 3.5. We define a permutation $\tau$ of
$\{ 1, \ldots , p+q \}$ to be in $\sanc$ when every cycle of $\tau$
is ``standard'' in the appropriate annular sense, and when $\tau$ does 
not display any of these crossing patterns. By starting from this definition, 
we then prove that $\sanc$ can also be described via an annular counterpart
of Biane's geodesic condition; this is done in the Theorem 6.1 of the paper.
On the way towards that, we observe that at the algebraic level there exists 
a nice connection between $\sanc$ and the set of disc non-crossing 
permutations $\GS_{disc-nc} (p+q)$. Namely: $\sanc$ is the saturation of 
$\GS_{disc-nc} (p+q)$ under conjugation with the permutations 
$\gamma_{ext}, \gamma_{int} \in \GS (p+q),$ where 
$\gamma_{ext} := (1, \ldots , p-1,p)$ and 
$\gamma_{int} := (p+1, \ldots , p+q-1,p+q)$; cf. Theorem 5.1.

\vspace{10pt}

It is well-known that disc non-crossing partitions (or permutations) play
a role in describing the asymptotics for moments in several important
examples of random matrices. The annular non-crossing permutations turn
out to play a similar role in the description of the second order
asymptotics for these matrices. We illustrate this phenomenon on a family of 
complex Wishart matrices (a family $G_1^* G_1, \ldots , G_s^* G_s$, where
$G_1, \ldots , G_s$ are $M \times N$ random matrices with independent 
complex Gaussian entries); this is done in the Section 7 of the paper. At
the end of the Section 7 we comment on how the same phenomenon also appears
in connection to Gaussian Hermitian random matrices; this example is not 
so illustrative for the purposes of the present paper, as it involves only 
the special case of complete matchings (permutations $\tau$ such that every
orbit of $\tau$ has exactly two elements), rather than dealing with 
general permutations.

\vspace{10pt}

In view of the analogy with the disc case, we cannot omit a discussion of 
the concept of annular non-crossing partition. This is done in the Section
4 of the paper. We define a partition $\pi$ of $\{ 1, \ldots , p+q \}$ to
be annular non-crossing if there exists $\tau \in \sanc$ such that $\pi$ is
the orbit partition of $\tau$. The set $\ncannpq$ of partitions of 
$\{ 1, \ldots , p+q \}$ which is obtained in this way coincides with the one 
studied by King \cite{Ki} (the difference between the approaches taken here 
and in \cite{Ki} is that our definition of $\ncannpq$ is algebraic, while
King's is phrased in topological terms -- cf. Remark 4.3.2).
We point out that the natural surjection from $\sanc$ onto $\ncannpq$ is 
not one-to-one, and we identify precisely the kind of pathology 
which causes this to happen.

Since (unlike what we had in the disc case) $\sanc$ and $\ncannpq$ cannot
be identified to each other, one has has to choose which kind of object
(partitions or permutations?) should play the primary role in the 
development of the annular non-crossing theory. In the support of 
our choice for permutations we can bring the following two arguments:

(a) The complementation map introduced by Kreweras in the disc case has
a very nice counterpart at the level of permutations (a bijection from 
$\sanc$ to itself -- see Remark 6.7). But when trying to define the 
Kreweras complementation as a map on $\ncannpq$, one finds (as noticed 
already in \cite{Ki}) that the map is not defined everywhere, and is 
not a bijection on the part of $\ncannpq$ where it is defined. This fact 
is directly related to the lack of injectivity of the natural map from
$\sanc$ to $\ncannpq$.

(b) The second order asymptotics for random matrices found in Section 7
involve summations naturally indexed by $\sanc$ (cf. Equation (7.5) 
in Theorem 7.5); these could be turned into summations over $\ncannpq$ only 
by introducing some artificial weights. 

\vspace{10pt}

In the Section 8 of the paper we briefly discuss some possible generalizations
of our results from the framework of the annulus to the more general one of 
a multi-annulus. Then in the final Section 9 we see how the formulas for the 
joint moments of a Wishart family extend to the multi-annular framework. As an
application, we observe how a phenomenon of asymptotic Gaussianity for traces 
of words made with independent Wishart matrices can be obtained via some 
fairly straightforward combinatorial arguments (cf. Corollary 9.4).

\vspace{10pt}

As pointed to us by Vaughan Jones, the combinatorics of annular non-crossing
permutations developped in this paper has striking resemblances with 
the combinatorics of ``tangles'', objects appearing in Jones' theory of 
planar algebras (see e.g. \cite{J2}). This coincidence suggests that the 
second-order asymptotics of random matrices could perhaps be of use in 
constructions of von Neumann subfactors (but at the present moment there
is no precise idea on how to substantiate such a possibility).

$\ $

$\ $

\section{Review of non-crossing permutations in the disc} 
 
We start from the better established concept of non-crossing partition.

$\ $

{\bf 2.1 Notations.} Let $n$ be a positive integer. We will denote 
the set $\{ 1, \ldots , n \}$ by $[n]$.

$1^o$ Let $\pi$ be a partition of $[n]$; that is,
$\pi = \{ B_1, \ldots , B_k \}$ where $B_1, \ldots , B_k$ (called the 
blocks of $\pi$) are non-empty subsets of $[n]$ such that 
$B_1 \cup \cdots \cup B_k = [n]$ and such that $B_i \cap B_j = \emptyset$
for $i \neq j$. We say that $\pi$ {\em has crossings} if there exist 
$1 \leq i,j \leq k,$ $i \neq j,$ and $1 \leq a<b<c<d \leq n,$ such that
$a,c \in B_i$ and $b,d \in B_j.$ We say that $\pi$ is {\em non-crossing} if
it has no crossings.

$2^o$ A partition $\pi = \{B_1, \ldots , B_k \}$ of $[n]$
can be represented pictorially on a disc in the following way: one first
draws the points $1, \ldots , n$ around the boundary of the disc (our choice 
is to do so clockwise), and then for every block $B$ of $\pi$ one draws the 
convex hull of $\{ i \ | \ i \in B \}.$ The partition $\pi$ is non-crossing 
precisely when the convex hulls drawn for its blocks are pairwise disjoint.

$3^o$ The set of non-crossing partitions of $[n]$ is usually denoted as 
$NC(n).$ In this paper we will denote it as $\ncdisc (n),$ in order to 
distinguish it from the set of non-crossing partitions in an annulus which 
will be discussed in Section 4.

$\ $

As explained in the Introduction, we will follow the approach of viewing
$\ncdisc (n)$ as embedded in the group of permutations of $[n]$. We next 
review some basic terminology related to this approach.

$\ $

{\bf 2.2 Notations.} $1^o$ For a finite non-empty set $A$, we will denote 
the group of all permutations of $A$ as $\GS (A).$

$2^o$ Let $A$ be a finite non-empty set and let $\tau$ be a permutation in
$\GS (A).$ Then $A$ is partitioned into {\em orbits} of $\tau$
($a_1 , a_2 \in A$ are in the same orbit of $\tau$ when there exists 
$k \in \Z$ such that $\tau^k (a_1) = a_2).$ We will denote
\begin{equation}
\# ( \tau ) \ : = \ \mbox{the number of orbits of $\tau$.}
\end{equation}
If $B \subset A$ is an orbit of $\tau$, then the restriction of $\tau$ to
$B$ (which is a permutation $\tau \mid B \in \GS (B)$) is called the
{\em cycle} of $\tau$ corresponding to the orbit $B$. A permutation 
$\tau \in \GS (A)$ is said to be {\em cyclic} if it has only one orbit,
equal to $A$.

$3^o$ We will usually write permutations in cycle notation. For instance,
if $A = \{ 1,3,4,5,6,$
$9 \}$ and if 
$\tau = \left( \begin{array}{cccccc} 1 & 3 & 4 & 5 & 6 & 9  \\
4 & 9 & 6 & 5 & 1 & 3  \end{array} \right)$, then we will write
``$\tau = (1,4,6)(3,9)(5)$'', or just 
``$\tau = (1,4,6)(3,9) \in \GS (A)$'' (as it is customary
to omit the orbits with one element from the cycle notation).

$\ $

{\bf 2.3 Remark and Definition} {\em (permutation induced  on a subset).} 
Let $A$ be a finite set, let $B$ be a non-empty subset of $A$, and 
let $\tau$ be a permutation in $\GS (A).$ The permutation 
``$\tau \mid B \in \GS (B)$'' can be defined naturally,
even if $B$ is not an orbit of $\tau$, or a union of such orbits. To be 
precise, $\sigma = \tau \mid B \in \GS (B)$ is defined in 
the following way: for every $b \in B$ we look at the sequence (in $A$)
$\tau (b), \tau^2 (b), \ldots , \tau^k (b), \ldots ,$ and define 
$\sigma (b)$ to be the first element of this sequence which is again in $B$.
In what follows we will refer to $\tau \mid B$ as ``{\em the permutation induced by $\tau$ on $B$}''.

It is immediate to check that (for $\tau \in \GS (A)$ and $B \subset A$ as 
above) the orbits of $\tau \mid B$ are precisely the sets of the form 
$C \cap B$, where $C$ is an orbit of $\tau$ which intersects $B$. In 
particular, $\tau \mid B$ is a cyclic permutation of $B$ if and only if 
$B$ is contained in one of the orbits of $\tau$.

Note also that the operation of inducing a permutation to a subset behaves
well when done repeatedly, i.e. we have 
$( \tau \mid B ) \mid C = \tau \mid C$
whenever $\emptyset \neq C \subset B \subset A$ and $\tau \in \GS (A).$

$\ $

{\bf 2.4 Definition and Remark} {\em (standard permutations in disc sense.)}
Let $n$ be a positive integer and let $\gamma_o$ be the forward cyclic 
permutation of $[n],$
\begin{equation}
\gamma_{o} \ := \ (1, \ldots , n-1 , n) \in \GSn .
\end{equation}
If a permutation $\tau \in \GSn$ has the property that 
\begin{equation}
\tau \mid B = \gamma_o \mid B, \mbox{ for every orbit $B$ of $\tau$,}
\end{equation}
then we will say that $\tau$ is {\em standard in disc sense.}
(The condition (2.3) asks, in other words, that for every orbit 
$B = \{ b_1, \ldots , b_k \}$ of $\tau$, with $b_1< \cdots < b_{k-1}<b_k$,
the corresponding cycle of $\tau$ is $(b_1, \ldots , b_{k-1}, b_k)$.)

$\ $

{\bf 2.5 Definition} {\em (embedding of $\ncdisc (n)$ into $\GSn$).} Let 
$n$ be a positive integer.

$1^o$ For every partition $\pi \in \ncdisc (n)$ we will denote as
$\perm_{\pi}$ the unique permutation of $[n]$ which is standard in disc 
sense and which has the blocks of $\pi$ as orbits.

$2^o$ The set $\{ \perm_{\pi} \ | \pi \in \ncdisc (n) \}$ will be denoted
as $\sdnc$, and its elements will be called
{\em disc non-crossing permutations.} 

$\ $

{\bf 2.6 Remarks.} $1^o$ Geometrically, the map
$\ncdisc (n) \ni \pi \mapsto \perm_{\pi} \in \sdnc$ is thus described as 
follows: in the disc picture of $\pi$ (as described in Definition 2.1.2) 
we ask that the boundary of each convex polygon drawn for a block of $\pi$ is 
run clockwise by the corresponding cycle of $\perm_{\pi}$. (See the 
Figure 1 at the end of the paper for a concrete example of such a drawing.)

$2^o$ On the other hand, let us pursue in more detail the algebraic 
description for the fact that a permutation $\tau$ of $[n]$ belongs to 
$\sdnc$. There are two things which are required:
(i) that $\tau$ is standard in disc sense, and
(ii) that the partition of $[n]$ into orbits of $\tau$ is in $\ncdisc (n).$
We would like to point out that the negations of both (i) and (ii) can be 
phrased in terms of some simple {\em localized} equations, involving $\tau$
and the forward cycle $\gamma_o$ of Equation (2.2). The term ``localized''
refers to the fact that the equations will only focus on what is induced by 
$\tau$ and $\gamma_o$ on subsets of $[n]$ with not more than 4 elements.
Before giving the precise statement of how this goes (in Proposition 2.8),
we record the following simple observation, the proof of which is left to 
the reader:

$\ $

{\bf 2.7 Lemma.} Let $B$ be a finite non-empty set, and let 
$\sigma_1 , \sigma_2$ be two cyclic permutations of $B$.
If $\sigma_1 \mid C = \sigma_2 \mid C$ for every subset $C \subset B$
which has 3 elements, then $\sigma_1 = \sigma_2.$

$\ $

{\bf 2.8 Proposition.} Let $n$ be a positive integer, and let $\gamma_o$ be
as in (2.2). Consider the following two conditions, which a permutation 
$\tau \in \GSn$ may or may not fulfill:

{\em (DNS)} There exist 3 distinct elements $a,b,c \in [n]$ such that 
$\gamma_o \mid \{ a,b,c \} = (a,b,c)$ and $\tau \mid \{ a,b,c \} = (a,c,b).$ 

{\em (DC)} There exist 4 distinct elements $a,b,c,d \in [n]$ such that 
$\gamma_o \mid \{ a,b,c,d \} = (a,b,c,d)$ and
$\tau \mid \{ a,b,c,d \} = (a,c)(b,d).$ 

Then for a permutation $\tau \in \GSn$ we have that: $\tau \not\in \sdnc$ 
if and only if $\tau$ satisfies at least one of the conditions (DNS) and (DC).

\vspace{6pt}

[All the restrictions of permutations appearing in (DNS) and in (DC) are 
in the ``induced'' sense discussed in Section 2.3. The acronyms ``DNS'' and 
``DC'' stand for ``Disc Non--Standard'' and respectively for ``Disc 
Crossing''.]

$\ $

{\bf Proof.} It is immediate from the definitions that $\tau$ satisfies 
(DC) if and only if the partition of $[n]$ into orbits of $\tau$ has 
crossings. On the other hand, when we use the Lemma 2.7 in connection to 
Equation (2.3), we see that $\tau$ satisfies (DNS) if and
only if it is not standard in disc sense. {\bf QED}

$\ $

{\bf 2.9 Remark.} Let $n$ be a positive integer, and let $T$ be a totally
ordered set with $n$ elements. The above discussion about non-crossing 
partitions/permutations of $[n]$ can be transferred to 
partitions/permutations of $T,$ by using the unique order-preserving bijection 
$\varphi : [n] \to T$. More precisely, we define:
\begin{center}
$\GS_{disc-nc} (T) \ = \ \{ \varphi \tau \varphi^{-1} \ | \ 
\tau \in \sdnc \} ,$

\vspace{6pt}

$NC_{disc} (T) \ = \ \Bigl\{ \ \{ \varphi (B_1), \ldots , \varphi (B_k) \} \ 
| \ \{ B_1, \ldots , B_k \} \in NC_{disc}(n) \ \Bigr\} .$
\end{center}
\noindent
It is immediate that the Proposition 2.8 holds for permutations in $\GS (T),$ 
if $\gamma_o$ is replaced with the forward cyclic permutation of $T$ (or in 
other words, with $\varphi \gamma_o \varphi^{-1}$).
 
$\ $

{\bf 2.10 Review} {\em of Biane's geodesic condition.} Let $n$ be a positive
integer. A remarkable fact observed by Biane (\cite{B2}, Section 1.3) is that
$\sdnc$ can also be described as follows: it is the set of permutations 
$\tau \in \GSn$ which satisfy the relation
\begin{equation}
\# ( \tau ) \, + \, \# ( \tau^{-1} \gamma_o ) \ = \ n+1,
\end{equation}
with $\# ( \tau )$ as defined in Equation (2.1), and 
$\gamma_o = (1, \ldots , n-1,n)$.
This is in the context where the inequality:
\begin{equation}
\# ( \tau ) \, + \, \# ( \tau^{-1} \gamma_o ) \leq n+1
\end{equation}
is satisfied by all the permutations $\tau \in \GSn$.

The Equation (2.4) can be viewed as a ``geodesic condition'', in the 
following sense. Consider the Cayley graph of the group $\GSn$, where the 
set of generators chosen for $\GSn$ is the set of all transpositions,
$\{ (i,j) \ | \ 1 \leq i<j \leq n \} \subset \GSn$. (This means that the 
set of vertices of the Cayley graph is $\GSn$, and that we draw an edge 
between the vertices $\sigma$ and $\tau$ when $\sigma^{-1} \tau$ is a 
transposition.) The distance in this Cayley graph turns out to be 
described by the formula
\begin{equation}
d( \sigma , \tau ) \ = \ n - \# ( \sigma^{-1} \tau ), \ \ 
\mbox{ for } \sigma , \tau \in \GSn .
\end{equation}
But then, as is immediately verified, the inequality (2.5) amounts to 
just the triangle inequality
\begin{equation}
d( \, id , \tau  \, ) + d( \, \tau , \gamma_o \, ) \ \geq \
d( \, id , \gamma_o \, ),
\end{equation}
where $id$ denotes the identity permutation of $[n].$ Hence in other words,
Biane's description of $\sdnc$ is thus: a permutation $\tau$ of $[n]$ is 
disc non-crossing if and only if it satisfies (2.7) with equality, i.e. if
and only if it lies on a geodesic connecting
$id$ and $\gamma_o$ in the Cayley graph of $\GSn$.

A slight generalization of the geodesic condition (also observed by Biane 
-- see \cite{B1}, Lemma 3) goes as follows: Suppose that $\omega \in \GSn$ 
is standard in disc sense, and let $B_1, \ldots , B_k$ be the orbits of 
$\omega$. Then for $\tau \in \GSn$ we have:
\begin{equation}
d( \, id , \tau \, ) + d( \, \tau , \omega \, ) = d( \, id , \omega \, )  \
\Leftrightarrow \ \left\{  \begin{array}{c}
B_j \mbox{ is a union of orbits of $\tau$, and }    \\
\tau \mid B_j \in \GS_{disc-nc} (B_j), \ \forall \ 1 \leq j \leq k
\end{array} \right\} .
\end{equation}
(In (2.8), $\GS_{disc-nc} (B_j)$ is considered in the sense of Remark 2.9,
where the total order on $B_j$ is the one induced from $[n].$)

$\ $

{\bf 2.11 Notation and Remark.} For $\tau , \sigma \in \GSn$ let us denote
\begin{equation}
\# ( \tau \vee \sigma ) \ := \ \begin{array}{l}
\mbox{the number of orbits into which $[n]$ is split}  \\
\mbox{by the joint action of $\tau$ and $\sigma$}
\end{array}
\end{equation}
(where by ``joint action of $\tau$ and $\sigma$'' we understand the action 
on $[n]$ of the subgroup of $\GSn$ generated by $\tau$ and $\sigma$). Note 
that we always have
$\# ( \tau \vee \sigma ) \ \leq \ \min ( \ \# ( \tau ), \# ( \sigma ) \ )$,
since a joint orbit of $\tau$ and $\sigma$ is a union of orbits of $\tau$, 
and also a union of orbits of $\sigma$. In the case when 
$\# ( \tau \vee \sigma ) = 1$ we will say that the joint action of 
$\tau$ and $\sigma$ is {\em transitive} on $[n]$.

A useful generalization of the inequality (2.5) is that:
\begin{equation}
\# ( \tau ) + \# ( \tau^{-1} \sigma ) + \# ( \sigma ) \ \leq \ 
n + 2 \cdot \# ( \tau \vee \sigma ), \ \ \forall \ \tau , \sigma \in \GSn .
\end{equation}
This is well-known, and appears in various forms in the literature
(see e.g. the Section 2 of \cite{GJ}). For the convenience of the reader, 
we briefly go over the steps of its proof. First, we use a summation over 
the set of joint orbits of $\tau$ and $\sigma$ in order to reduce (2.10) to 
its particular case when $\tau$ and $\sigma$ act transitively on $[n]$. 
Then in the case when $\tau$ and $\sigma$ act transitively we proceed by 
induction on $\# ( \sigma )$. The base case $\# ( \sigma ) =1$ is essentially
the inequality (2.5). For the induction step: suppose that the case when 
$\# ( \sigma ) = m-1$ was settled, and that we want to prove the case when 
$\# ( \sigma ) = m$ (for some $2 \leq m \leq n$). So let $\sigma , \tau 
\in \GSn$ be such that $\# ( \sigma ) = m$ and such that the joint action
of $\tau$ and $\sigma$ on $[n]$ is transitive. There have to exist 
$1 \leq i < j \leq n$ such that $i$ and $j$ belong to the same orbit
of $\tau$, but are in distinct orbits of $\sigma$. (Indeed, in the opposite
case it would follow that every orbit of $\tau$ is contained in an orbit of 
$\sigma$, and this would immediately contradict the transitivity of the 
joint action of $\tau$ and $\sigma$.) Consider the permutations 
$\widetilde{\sigma} := (i,j) \sigma$ and $\widetilde{\tau} := (i,j) \tau$.
It is immediately verified that 
$\# ( \widetilde{\sigma} ) = \# ( \sigma ) -1 = m-1$, that
$\# ( \widetilde{\tau} ) = \# ( \tau ) +1$, and that the joint action
of $\widetilde{\sigma}$ and $\widetilde{\tau}$ on $[n]$ is transitive. As a 
consequence we get that
\[
\# ( \tau ) + \# ( \tau^{-1} \sigma ) + \# ( \sigma ) \ = \ 
\# ( \widetilde{\tau} ) + \# ( \widetilde{\tau}^{-1} \widetilde{\sigma} )
+ \# ( \widetilde{\sigma} ) \ \leq \  n+2,
\]
as desired, where the last inequality holds by the induction hypothesis. 

$\ $

$\ $

\section{Annular non-crossing permutations via localized conditions} 

In this section we fix two positive integers 
\setcounter{equation}{0}
$p$ and $q$, and we will look at permutations in $\GSpq$.

As the name suggests, the concept of $(p,q)$--annular non-crossing
permutation will be defined by identifying what are the crossing patterns 
in $(p,q)$--annular sense; then a permutation $\tau$ of $[p+q]$ will be 
declared to be $(p,q)$--annular non-crossing if it is ``standard'' (in the 
appropriate annular sense) and if it does not display any of these crossing
patterns. The definition will thus be algebraic, in terms of ``localized 
conditions'' (if $\tau \in \GSpq$ is not annular non-crossing, then this will 
be detectable by inspecting a group of not more than 6 elements of $[p+q],$ 
which belong to not more than 3 distinct orbits of $\tau$, and by checking for 
these elements some conditions analogous to (DNS) and (DC) of Proposition 2.8).

In this section we will also explain (at least on a heuristic level) why 
our algebraic definition does indeed lead to the set of permutations which 
one expects to see when drawing pictures. In order to give a better intuition
of what is going on, we start by explaining what kind of ``pictures of
permutations'' we have in mind. In discussions based on pictures we will use 
the term ``planar ''(and we will reserve the term ``non-crossing'' for the 
algebraic approach based on crossing patterns).

$\ $

{\bf 3.1 Remark} {\em (planarity in the annular framework).} Instead of 
visualising permutations by drawing $n = p+q$ points around a circle, we will
now use two concentric circles. The external circle has marked on it the
points $\{ 1, \ldots , p \}$, in clockwise order (matching the choice of 
running clockwise made in Section 2). The internal circle has marked on 
it the points $\{ p+1, \ldots , p+q \}$, in counter-clockwise order.
Our drawings will be made in the annulus between the two circles, which 
will be referred to as ``{\em the $(p,q)$--annulus}''. (Note: marking the 
points $\{ p+1, \ldots , p+q \}$ counter-clockwise is the consistent 
choice for someone who lives in the $(p,q)$--annulus.) 

Let $\sigma = (a_1, \ldots, a_k)$ be a cyclic permutation of a subset
$A \subset [p+q]$. By a {\em picture of $\sigma$ in the $(p,q)$--annulus}
we will understand a curve in the $(p,q)$--annulus which connects $a_1$ to
$a_2$, then $a_2$ to $a_3, \ldots ,$ then $a_k$ to $a_1$, in such a way that
it only touches the boundary of the annulus at the points $a_1, \ldots , a_k$.
Such a picture will be said to be {\em admissible} if it does not 
self-intersect, if it encloses a region which is contained in the 
$(p,q)$--annulus, and if it winds clockwise around that region.

Now let $\tau$ be a permutation of $[p+q]$. If we can draw an admissible 
picture for every cycle of $\tau$ (in the sense of the preceding paragraph)
such that the regions enclosed by two different cycles are disjoint, then 
we will say that we have a {\em planar $(p,q)$--annular picture} of $\tau$.
A concrete example of such a drawing is shown in the Figure 2 at the end 
of the paper.

The use of the name ``planar'' is justified by the fact that (at least in 
the most interesting case, when $\tau$ has at least one orbit which 
intersects both $\{ 1, \ldots , p \}$ and $\{ p+1, \ldots , p+q \}$) 
our drawing is a 
planar map in the sense of Tutte (see e.g. \cite{Tu2}). It is convenient to
imagine that the faces of the map are tri-coloured in black, white and red:
There are 2 red faces (namely the interior of the internal circle and the 
exterior of the external circle of the $(p,q)$--annulus), and there is one
black face for every orbit of $\tau$ (the checked areas in Figure 2); then
the white faces cover what is left of the $(p,q)$--annulus.

$\ $

We now start on the algebraic approach, via crossing patterns. The algebraic
counterparts for the topological terms appearing in Remark 3.1 will be as
follows:

\vspace{6pt}

\begin{center}
\begin{tabular}{c|c}   \hline
                    &                       \\
Topological terms   &  Algebraic terms      \\
                    &                       \\  \hline
                    &                       \\
Every cycle of $\tau$ has an admissible & 
$\tau$ is standard in $(p,q)$--annular sense  \\
picture in the $(p,q)$--annulus         &
 (cf. Definition 3.3, Remarks 3.4)           \\
                    &                       \\  \hline
                    &                       \\
$\tau$ has a planar $(p,q)$--annular picture & 
$\tau$ is $(p,q)$--annular non-crossing (cf. Def. 3.5)  \\ 
                    &                       \\  \hline
\end{tabular}
\end{center}

\vspace{6pt}

The explanation of why the correspondence is as stated in the above table
will be given in the Sections 3.6-3.9 below.

$\ $

{\bf 3.2 Notations.} $1^o$ We will extensively use the permutations 
$\gamma_{ext}, \gamma_{int}, \gamma \in \GSpq$ defined by:
\begin{equation}
\left\{  \begin{array}{ll}
\gamma_{ext} := (1, \ldots , p-1,p),  &
\gamma_{int} := (p+1, \ldots , p+q-1 ,p+q),   \\
                                     &         \\
\gamma  := \gamma_{ext} \gamma_{int} = \gamma_{int} \gamma_{ext} &
   = (1, \ldots , p-1,p) \, (p+1, \ldots , p+q-1, p+q).
\end{array}  \right.
\end{equation}

$2^o$ Let $x$ be in $\{ 1, \ldots , p \}$ and let $y$ be in 
$\{ p+1, \ldots , p+q \}.$ Consider the permutation in $\GSpq$ which fixes
$x$ and $y$, and which organizes $[p+q] \setminus \{ x,y \}$ in a cycle as
follows:
\begin{equation}
\begin{array}{ccc}
\gamma_{ext}(x)   & \rightarrow   \gamma_{ext}^2(x) \rightarrow \cdots
                                \rightarrow      &  \gamma_{ext}^{p-1}(x)  \\
                       &                        &                         \\
\uparrow               &                        &  \downarrow             \\
                       &                        &                         \\
\gamma_{int}^{q-1}(y) & \leftarrow \cdots \leftarrow \gamma_{int}^2 (y)
                               \leftarrow      &  \gamma_{int} (y)
\end{array}
\end{equation}
This permutation will be denoted as $\lambda_{x,y}$ (or sometimes as 
$\lambda_{y,x}$, if so needed).

The intuitive significance of $\lambda_{x,y}$ is thus: if we cut the 
$(p,q)$--annulus along a simple curve which connects $x$ and $y$, then
we obtain a topological disc, on the boundary of which the points from 
$[p+q] \setminus \{ x,y \}$ sit in the order indicated by $\lambda_{x,y}$
(see Figure 3).

On the algebraic side, note that one can introduce $\lambda_{x,y}$ by more 
concise formulas than (3.2); for instance it is immediately verified that,
for $p,q \geq 2,$ one can also write
\begin{equation}
\lambda_{x,y} \ = \ \gamma_{int} \cdot
( \, \gamma_{ext}(x), x, y , \gamma_{int}^{-1} (y) \, ) \cdot \gamma_{ext}.
\end{equation}

$\ $

{\bf 3.3 Definition.} We will say that a permutation $\tau \in \GSpq$ is
{\em standard in $(p,q)$--annular sense} when it satisfies the following
conditions:

(i) For every orbit $A$ of $\tau$ we have that 
$\tau \mid A \cap \{ 1, \ldots , p \}$
= $\gamma_{ext} \mid A \cap \{ 1, \ldots , p \}$ and that 
$\tau \mid A \cap \{ p+1, \ldots , p+q \}$ = 
$\gamma_{int} \mid A \cap \{ p+1, \ldots , p+q \}.$

(ii) If $A$ is an orbit of $\tau$ such that $A \cap \{ 1, \ldots , p \}
\neq \emptyset \neq A \cap \{ p+1, \ldots , p+q \}$, then there exist 
a unique $a' \in A \cap \{ 1, \ldots , p \}$ such that 
$\tau (a') \in \{ p+1, \ldots , p+q \}$ and a unique 
$a'' \in A \cap \{ p+1, \ldots , p+q \}$ such that 
$\tau (a'') \in \{ 1, \ldots , p \}$.

$\ $

{\bf 3.4 Remarks.} $1^o$ Let $\tau \in \GSpq$ be standard in 
$(p,q)$--annular sense and let $A$ be an orbit of $\tau$. Then the 
possibilities for $\tau \mid A$ are quite limited.
We distinguish two cases:

\vspace{4pt}

{\em Case 1:} $A \subset \{ 1, \ldots , p \}$, or 
$A \subset \{ p+1 , \ldots , p+q \}$.

In this case we have that $\tau \mid A = \gamma \mid A$ (where the latter
permutation is equal in turn to $\gamma_{ext} \mid A$ if 
$A \subset \{ 1, \ldots , p \}$, and to $\gamma_{int} \mid A$ if
$A \subset \{ p+1, \ldots , p+q \}$).

\vspace{4pt}

{\em Case 2:} $A \cap \{ 1, \ldots , p \} \neq \emptyset \neq 
A \cap \{ p+1 , \ldots , p+q \}$.

In this case we have that
\begin{equation}
\tau \mid A \ = \ ( a_1, \ldots , a_k, b_1, \ldots , b_l ),
\end{equation}
where $a_1, \ldots , a_k \in \{ 1, \ldots , p \}$ are such that 
$\gamma_{ext} \mid \{ a_1, \ldots , a_k \} = ( a_1, \ldots , a_k),$ and
where $b_1, \ldots , b_l \in \{ p+1, \ldots , p+q \}$ are such that 
$\gamma_{int} \mid \{ b_1, \ldots , b_l \} = ( b_1, \ldots , b_l).$ 
(Note: the latter conditions are saying that each of the 
sequences $a_1, \ldots , a_k$ and $b_1, \ldots , b_l$, taken separately,
can be permuted cyclically to be put in increasing order. This, of course,
isn't generally true for the whole sequence of $a$'s and $b$'s in (3.4) 
-- it could for instance happen that 
$\tau \mid A = (p-1,p,1, p+q-1, p+q, p+1)$.) A typical cycle of the kind 
appearing in (3.4) is drawn in the Figure 4 at the end of the paper.

\vspace{6pt}

$2^o$ We leave it as an exercise to the reader to check that a permutation
$\tau \in \GSpq$ is {\em not} standard in $(p,q)$--annular sense if and only 
if it satisfies at least one of the following two ``localized'' conditions:

\vspace{4pt}

{\em (ANS-1)} There exist 3 distinct elements $a,b,c \in [p+q]$ such that 
$\gamma \mid \{ a,b,c \} = (a,b,c)$ and $\tau \mid \{ a,b,c \} = (a,c,b).$ 

\vspace{4pt}

{\em (ANS-2)} There exist 4 distinct elements $a,b,c,d \in [p+q]$ such that 
$\gamma \mid \{ a,b,c,d \} = (a,b)(c,d)$ and 
$\tau \mid \{ a,b,c,d \} = (a,c,b,d).$ 

$\ $

{\bf 3.5 Definition.} Consider the following conditions, which a permutation
$\tau \in \GSpq$ may or may not fulfill.

{\em (AC-1)} There exist 4 distinct elements $a,b,c,d \in [p+q]$ such that 
$\gamma \mid \{ a,b,c,d \} = (a,b,c,d)$ and
$\tau \mid \{ a,b,c,d \} = (a,c)(b,d).$ 

{\em (AC-2)} There exist 5 distinct elements $a,b,c,x,y \in [p+q]$ such that 
$x$ and $y$ belong to the two distinct orbits of $\gamma$, and such that 
$\lambda_{x,y} \mid \{ a,b,c \} = (a,b,c)$ and
$\tau \mid \{ a,b,c,x,y \} = (a,c,b)(x,y).$ 

{\em (AC-3)} There exist 6 distinct elements $a,b,c,d,x,y \in [p+q]$ such 
that $x$ and $y$ belong to the two distinct orbits of $\gamma$, and such 
that $\lambda_{x,y} \mid \{ a,b,c,d \} = (a,b,c,d)$ and
$\tau \mid \{ a,b,c,d,x,y \} = (a,c)(b,d)(x,y).$ 

A permutation $\tau \in \GSpq$ will be said to be
{\em $(p,q)$--annular non-crossing} if it is standard in $(p,q)$--annular 
sense, and if it does not satisfy any of the conditions (AC-1), (AC-2), 
(AC-3). The set of all such permutations will be denoted as $\sanc$.

$\ $

{\bf 3.6 Remark.} Let $\tau$ be a permutation of $[p+q].$ As announced at 
the beginning of the section, the statement ``$\tau \not\in \sanc$'' amounts
to the fact that $\tau$ satisfies at least one of the 5 localized conditions 
(ANS-1), (ANS-2), (AC-1), (AC-2), (AC-3). The acronyms ``ANS''
and ``AC'' stand for ``Annular Non-Standard'' and respectively for 
``Annular Crossing''.

So (AC-1), (AC-2), (AC-3) are the three crossing patterns which can be 
exhibited by a permutation $\tau$ of $[p+q]$ such that $\tau$ is standard
in $(p,q)$--annular sense, but $\tau \not\in \sanc$. (Note that these three
crossing patterns aren't exclusive to each other; for instance if $a,b,c,d$
are as in (AC-1) and if we can find $x \in \{1, \ldots , p \},$
$y \in \{ p+1, \ldots , p+q \}$ such that $x,y$ belong to the same orbit 
$A$ of $\tau$, with $A \subset [p+q] \setminus \{ a,b,c,d \}$, then 
$a,b,c,d,x,y$ will be as in (AC-3).) The crossing patterns (AC-1), (AC-2), 
(AC-3) are illustrated in the Figures 5--7 at the end of the paper. In the 
rest of the remark we make a brief comment on each of them.

The crossing pattern (AC-1) is very similar to the pattern (DC) from the disc
case. (Note that if $\gamma \mid \{ a,b,c,d \} = (a,b,c,d)$, then all 4 of 
$a,b,c,d$ belong to the same orbit of $\gamma$, so we are actually dealing 
either with $\gamma_{ext} \mid \{ a,b,c,d \}$ or 
with $\gamma_{int} \mid \{ a,b,c,d \}$.)

When discussing the crossing pattern (AC-2), it is useful to make the 
following observation:
\begin{equation}
\left\{  \begin{array}{l}
\mbox{ If $\tau \in \GSpq$ is standard in $(p,q)$--annular sense, }   \\
\mbox{ and satisfies (AC-2) for $a,b,c,x,y \in [p+q],$ }              \\
\mbox{ then $a,b,c$ cannot all belong to the same orbit of $\gamma$. }
\end{array} \right.
\end{equation} 
Indeed, $\gamma$ and $\lambda_{x,y}$ induce the same permutation on 
$\{ 1, \ldots , p \} \setminus \{ x,y \}$ and on 
$\{ p+1, \ldots , p+q \} \setminus \{ x,y \}$; so if $a,b,c$ were in the 
same orbit of $\gamma$, then it would follow that
$\gamma \mid \{ a,b,c \} = \lambda_{x,y} \mid \{ a,b,c \} = (a,b,c),$ 
and $\tau$ would not be standard in $(p,q)$--annular sense (since 
$\tau \mid \{ a,b,c \} = (a,c,b)$). In view of (3.5), the drawings 
illustrating (AC-2) are as in Figure 6.

It is instructive to compare the Figure 6 with the Figure 6'. The
condition (AC-2) does not apply to $a,b,c,x,y$ of Figure 6', since $\tau$
and $\lambda_{x,y}$ induce the same permutation on $\{ a,b,c \}$. Pictorially,
this corresponds to the fact that in Figure 6' we have no crossings (where
the cycles appearing in both the Figures 6 and 6' are drawn so that 
they wind clockwise around the regions they enclose).

Finally,  the meaning of (AC-3) is thus: suppose that $\tau$ satisfies 
(AC-3) for $a,b,c,d,x,y$, and that we could draw the cycles of $\tau$
which contain $a,b,c,d,x,y$ such that neither the cycle containing $a,c$ 
nor the cycle containing $b,d$ crosses the cycle
containing $x,y$. Then we can cut the $(p,q)$--annulus along a simple 
curve which connects $x$ and $y$, and which does not intersect the cycles
of $\tau$ containing $a,c$ and $b,d$, respectively. What results is a 
topological disc with the points from $[p+q] \setminus \{ x,y \}$ 
distributed cyclically on its boundary, in the order indicated by 
$\lambda_{x,y}.$ But then the equation
$\lambda_{x,y} \mid \{ a,b,c,d \} = (a,b,c,d)$ (which is part of (AC-3))
shows that there is a crossing between the cycle of $\tau$ containing
$a,c$, and the cycle of $\tau$ containing $b,d$ (see Figure 7, where 
we illustrated the situations not falling under the incidence of (AC-1)).

$\ $

Now, let $\tau$ be a permutation of $[p+q]$ which is standard in 
$(p,q)$-annular sense. The preceding remark argues that if $\tau$ satisfies
at least one of the conditions (AC-1), (AC-2), (AC-3), then one cannot draw
a planar $(p,q)$--annular picture of $\tau$ (in the sense of
Remark 3.1). In order to complete the defense of Definition 3.5 we must also
consider the opposite situation, when $\tau$ does not satisfy any of (AC-1), 
(AC-2), (AC-3) (hence when $\tau \in \sanc$), and argue, at least 
heuristically, that in this situation $\tau$ {\em has} a planar
$(p,q)$--annular picture. It is convenient to do this separately in the 
cases when $\tau$ is (respectively is not) $(p,q)$--connected, in the
sense of the following definition.

$\ $

{\bf 3.7 Definition.}
Let $\tau$ be a permutation in $\GSpq$. We will say that an orbit 
$B$ of $\tau$ is {\em $(p,q)$--connecting} if 
$B \cap \{ 1, \ldots , p \} \neq \emptyset$ and
$B \cap \{ p+1, \ldots , p+q \} \neq \emptyset$. We will say that $\tau$ is
{\em $(p,q)$--connected} if it has at least one $(p,q)$--connecting orbit.

$\ $

{\bf 3.8 Remark.} Consider a permutation $\tau$ of $[p+q]$ which is
$(p,q)$--disconnected, and denote
$\tau \mid \{ 1, \ldots , p \} =: \tau_{ext}$ and
$\tau \mid \{ p+1, \ldots , p+q \} =: \tau_{int}$. Then
\begin{equation}
\tau \in \sanc \ \Leftrightarrow \ \left\{  \begin{array}{l}
\tau_{ext} \in \GS_{disc-nc} ( \, \{ 1, \ldots , p \} \, ) \mbox{ and}  \\
\tau_{int} \in \GS_{disc-nc} ( \, \{ p+1, \ldots , p+q \} \, ).
\end{array}  \right.
\end{equation}
Indeed, from the discussion in Remark 3.4.1 (Case 1) it is clear that $\tau$
is standard in $(p,q)$--annular sense if and only if both $\tau_{ext}$ and 
$\tau_{int}$ are standard in disc sense. Supposing that this happens,
we next observe that (AC-2)and (AC-3) don't apply at all to $\tau$, while
(AC-1) splits into two separate conditions on $\tau_{ext}$ and $\tau_{int}$, 
leading precisely to (3.6).

As a consequence of (3.6), we note that for a permutation of $[p+q]$ which
is $(p,q)$--disconnected we have: 
$\tau \in \sanc \Leftrightarrow \tau \in \GS_{disc-nc}(p+q).$

If $\tau \in \sanc$ is $(p,q)$--disconnected, then one draws a planar 
$(p,q)$--annular picture of $\tau$ by drawing planar pictures of 
$\tau_{ext}$ and of $\tau_{int}$ which stay very close to the external and 
respectively the internal circle of the $(p,q)$--annulus. (See Figure 8 at
the end of the paper.)

$\ $

{\bf 3.9 Remark.} Let $\tau \in \sanc$ be $(p,q)$--connected, and let us fix
a $(p,q)$--connecting orbit $A$ of $\tau$. The cycle $\tau \mid A$ is written 
explicitly as in (3.4), $\tau \mid A$ = 
$( a_1, \ldots , a_k,b_1, \ldots , b_l);$ the admissible picture of 
$\tau \mid A$ is as in Figure 4 at the end of the paper. Note that 
when the area enclosed by $\tau \mid A$ (checked area in Figure 4) is
shrunk to a point, its complement with respect to the $(p,q)$--annulus
is turned into a bouquet of $k+l-1$ discs. We will argue that:

(j) Every orbit $B \neq A$ of $\tau$ must be completely contained in one
of the $k+l-1$ discs of the bouquet. 

(jj) Suppose that $B$ and $B'$ are two distinct orbits of $\tau$, both 
distinct from $A$, and which are contained in the same disc of the above 
mentioned bouquet. Then $B$ and $B'$ are non-crossing inside that disc.

At least on a ``picture-based'' level, the assertions (j)$+$(jj) explain why 
the regions enclosed by the admissible pictures of the cycles of $\tau$ are
pairwise disjoint (and thus give together a planar $(p,q)$--annular
picture of $\tau$).

It remains to explain why (j) and (jj) do indeed hold.
To this end, we draw again the Figure 4 which represents $\tau \mid A,$
and this time we mark on it how the set $[p+q] \setminus A$ is split into 
$k+l-1$ subsets (some of them possibly empty) sitting inside the $k+l-1$
discs of the bouquet. We denote these $k+l-1$ subsets of $[p+q] \setminus A$
as $E_1, \ldots , E_{k-1}, I_{1}, \ldots , I_{l-1}, E_k \cup I_l$ -- see
Figure 9.

Suppose that $B$ is an orbit of $\tau$, $B \neq A,$ and let $x,y$ be two 
distinct points of $B$. Then the following things cannot happen:

$( \alpha_1 )$ $x \in E_i$ and $y \in E_j$ for $1 \leq i,j \leq k,$
$i \neq j.$ Indeed, in this case we could find two elements 
$a,a' \in A \cap \{ 1, \ldots ,p \}$ such that $\tau$ satisfies (AC-1) for
$a,a',x,y$.

$( \alpha_2 )$ $x \in I_i$ and $y \in I_j$ for $1 \leq i,j \leq l,$
$i \neq j.$ Same argument as for $( \alpha_1 ).$

$( \alpha_3 )$ One of $x,y$ is in $E_i$ for some $1 \leq i \leq k-1,$ while 
the other is in $I_1 \cup \cdots \cup I_{l}$
\newline
$\Bigl( \, = \{ p+1, \ldots , p+q \} \setminus A \, \Bigr).$ 
Indeed, in this case $\tau$ would satisfy (AC-2) for the 
5 elements $a_i, a_{i+1},b_1, x,y$.

$( \alpha_4 )$ One of $x,y$ is in $I_i$ for some $1 \leq i \leq l-1,$ while 
the other is in $E_1 \cup \cdots \cup E_k$. Same argument as for 
$(\alpha_3 ).$

After ruling out the possibilities described by $( \alpha_1 ) - ( \alpha_4 ),$
one only remains with the possibilities that $x,y \in E_i$ for some 
$1 \leq i \leq k-1,$ or that $x,y \in I_j$ for some $1 \leq j \leq l-1,$ or
that $x,y \in E_k \cup I_l.$ This verifies the above assertion (j).

Finally, let $B,B'$ be as in the assertion (jj). If $B,B' \subset E_i$
for some $1 \leq i \leq k-1,$ or if $B,B' \subset I_j$ for some 
$1 \leq j \leq l-1,$ then the fact that $B$ and $B'$ do not cross follows
 from the fact that $\tau$ does not fulfill (AC-1). 
If $B,B' \subset E_k \cup I_l,$ then the fact that $B$ and $B'$ do not cross
follows from the non-fulfillment of (AC-3), used with $a,c \in B,$ 
$b,d \in B',$ and with $x := a_k,$ $y := b_1$ taken from $A$.

$\ $

$\ $

\section{Annular non-crossing partitions} 

In order to respect the parallelism with the discussion made in the
\setcounter{equation}{0}
disc case, we will now briefly address the concept of annular non-crossing 
partition.

$\ $

{\bf 4.1 Definitions.} $1^o$ In this section we keep fixed the same positive 
integers $p,q$ as in Section 3. We will adopt the terminology introduced
in Section 3 concerning ``the $(p,q)$--annulus'' and the set $\sanc$ of 
$(p,q)$--annular non-crossing permutations.

$2^o$ If $\tau$ is a permutation of $[p+q]$, we will denote as 
``$\orbits ( \tau )$'' the partition of $[p+q]$ into orbits of $\tau$.

$3^o$ Let $\pi$ be a partition of $[p+q].$ We say that a block $B$ of $\pi$ is
{\em $(p,q)$--connecting} if $B \cap \{ 1, \ldots , p \} \neq \emptyset$ and
$B \cap \{ p+1, \ldots , p+q \} \neq \emptyset$. We say that $\pi$ is 
{\em $(p,q)$--connected} if it has at least one $(p,q)$--connecting block.

$\ $

It is clear that a permutation $\tau \in \GSn$ is $(p,q)$--connected (in the 
sense of Definition 3.7) if and only if the partition $\orbits ( \tau )$ 
is so.

$\ $

{\bf 4.2 Definition.} Let $\pi$ be a partition of $[p+q].$ We will say that 
$\pi$ is {\em $(p,q)$--annular non-crossing} if there exists $\tau \in \sanc$ 
such that $\pi = \orbits ( \tau ).$ The set of $(p,q)$-annular non-crossing 
partitions of $[p+q]$ will be denoted as $\ncannpq$.

$\ $

{\bf 4.3 Remarks.} $1^o$ The algebraic definition of $\ncannpq$ goes thus 
as follows: a partition $\pi$ of $[p+q]$ is in $\ncannpq$ if and only if 
it is possible that for every block $B$ of $\pi$ we choose a cyclic 
permutation of $B$ which is as in either Case 1 or Case 2 of Remark 3.4.1, 
in such a way that the resulting permutation of $[p+q]$ does not display 
any of the crossing patterns (AC-1), (AC-2), (AC-3).
It is conceivable that the algebraic definition of $\ncannpq$ can be 
simplified, but we are not aware of a nice way of doing so, at the present
moment. (An immediate difficulty is that the crossing pattern (AC-2) does
not seem to be transferable from permutations to partitions.)

$2^o$ At the looser level of pictorial definitions, we believe that 
$\ncannpq$ introduced in Definition 4.2 coincides with the set of partitions 
described by King in \cite{Ki}, and denoted there as ``$NC_2^A (p,q)$''. The 
definition of $NC_2^A (p,q)$ is formulated as follows: a partition $\pi$ of
$[p+q]$ is in $NC_2^A (p,q)$ when one can draw some arcs in the 
$(p,q)$--annulus, which connect the points marked $1, \ldots , p+q,$ such
that: (a) each block of $\pi$ is path-connected, and (b) arcs for different
blocks do not intersect. (See \cite{Ki}, Section 4.1, p. 1078.)

The inclusion $\ncannpq \subset NC_2^A (p,q)$ can be argued as follows:
Suppose that $\pi \in \ncannpq$, and let $\tau \in \sanc$ be
such that $\pi = \orbits ( \tau )$. Consider a planar $(p,q)$--annular 
picture of $\tau$ (as discussed in Section 3.1). Then by following 
the contours of the cycles of $\tau$ we obtain a family of arcs as required 
in the preceding paragraph, and thus we obtain that $\pi \in NC_2^A (p,q).$

For the opposite inclusion, $NC_2^A (p,q) \subset \ncannpq$, one needs to 
show that the family of arcs drawn for a partition $\pi \in NC_2^A (p,q)$ 
can always be chosen such that they give a planar $(p,q)$--annular picture
of a permutation $\tau \in \sanc$. 
A way to argue that this can be done is by simply ``thickening'' the arcs 
drawn for a $\pi \in NC_2^A (p,q)$, so that they can be identified with 
some contours which enclose very small areas. 
(The fact that a partition $\pi \in NC_2^A (p,q)$ can be drawn in this 
special way is implicitly appearing in the appendix of \cite{Ki} --  cf. the 
section of the appendix called ``Making $k$-bridges'', pp. 1087-1088 in 
\cite{Ki}.) 

$\ $

In the remaining of this section we will only use the algebraic definition
of $\ncannpq$. Our main observation is the following.

$\ $

{\bf 4.4 Proposition.} Let $\pi$ be a partition in $\ncannpq$, and suppose 
that $\pi$ has at least two $(p,q)$-connecting blocks. Then there exists a
{\em unique} permutation $\tau \in \sanc$ such that $\orbits ( \tau ) = \pi$.

\vspace{10pt}

{\bf Proof.} Let $\tau_1 , \tau_2 \in \sanc$ be such that 
$\orbits ( \tau_1 ) = \pi = \orbits ( \tau_2 ).$ In order to show that 
$\tau_1 = \tau_2$, it suffices that we fix a block $B$ of $\pi$, and show 
that $\tau_1 \mid B = \tau_2 \mid B$. If $B$ has less than 
3 elements, then the desired equality is clear (because $B$ has a unique 
cyclic permutation); so we will assume that $B$ has at least 3 elements.

Since $\pi$ has at least two $(p,q)$--connecting blocks, we can find such
a block $A$ which is distinct from $B$. We choose elements 
$x \in A \cap \{ 1, \ldots , p \}$ and $y \in A \cap \{ p+1, \ldots , p+q \}$,
and we look at the cycle $\lambda_{x,y}$ defined as in the Notations 3.2.

Let us observe that for any 3 distinct elements $a,b,c \in B$ we must have 
$\tau_1 \mid \{ a,b,c \} = \lambda_{x,y} \mid \{ a,b,c \}$ -- otherwise 
$\tau_1$ would satisfy the condition (AC-2) for the elements $a,b,c,x,y$.
Hence $\tau_1 \mid B$ and $\lambda_{x,y} \mid B$ induce the same permutation
on every 3-element subset of $B$, and the Lemma 2.7 gives us that 
$\tau_1 \mid B = \lambda_{x,y} \mid B.$ 

A similar argument shows that $\tau_2 \mid B = \lambda_{x,y} \mid B,$ and
this concludes the proof. {\bf QED}

$\ $

We next examine what are the counterparts of Proposition 4.4 in the cases of
partitions $\pi \in \ncannpq$ which have either 0 or 1 $(p,q)$--connecting 
blocks.

$\ $

{\bf 4.5 Proposition.} Let $\pi$ be a partition of $[p+q]$ which is
$(p,q)$-disconnected. Then $\pi$ is of the form 
$\pi = \pi_{ext} \cup \pi_{int},$
with $\pi_{ext}$ a partition of $\{ 1, \ldots , p \}$ and 
$\pi_{int}$ a partition of $\{ p+1, \ldots , p+q \}$. We have that:

$1^o$ $\pi \in \ncannpq$ if and only if 
$\pi_{ext} \in \ncdisc ( \, \{ 1, \ldots , p \} )$ and
$\pi_{int} \in \ncdisc ( \, \{ p+1, \ldots , p+q \} )$.

$2^o$ If $\pi \in \ncannpq$, then there exists a unique permutation 
$\tau \in \sanc$ such that $\orbits ( \tau ) = \pi$. This $\tau$ is obtained
by putting together the unique non-crossing permutations 
$\tau_{ext} \in \GS_{disc-nc} ( \, \{ 1, \ldots , p \} \, )$ and 
$\tau_{int} \in \GS_{disc-nc} ( \, \{ p+1, \ldots , p+q \} \, )$ which 
have $\orbits ( \tau_{ext} ) = \pi_{ext}$ and respectively
$\orbits ( \tau_{int} ) = \pi_{int}$.

\vspace{10pt}

{\bf Proof.} This is an immediate consequence of the description of the 
permutations in $\sanc$ which are $(p,q)$--disconnected, as discussed in 
Remark 3.8. The uniqueness of $\tau$ in $2^o$ follows for instance from the 
fact that every block of $\pi$ has a unique cyclic permutation which 
is standard in $(p,q)$--annular sense (cf. Remark 3.4.1, Case 1). {\bf QED}

$\ $

{\bf 4.6 Proposition.} Let $\pi$ be a partition of $[p+q]$ which has exactly
one $(p,q)$-connecting block $B_0$. We denote 
$B_0 \cap \{ 1, \ldots , p \} =: B_0 '$ and
$B_0 \cap \{ p+1, \ldots , p+q \} =: B_0 ''$. Also, let us denote the blocks
of $\pi$ which are completely contained in $\{ 1, \ldots , p \}$ 
(respectively in $\{ p+1, \ldots , p+q \}$) as $B_1', \ldots , B_s'$
(respectively as $B_1'', \ldots , B_t''$). We have that:

$1^o$ The number of permutations $\tau \in \GSpq$ which are standard in 
$(p,q)$--annular sense and have $\orbits ( \tau ) = \pi$ is equal to 
$\card ( B_0') \cdot \card ( B_0'').$

$2^o$ $\pi$ is in $\ncannpq$ if and only if the partition 
$\{ B_0', B_1', \ldots , B_s' \}$ of $\{ 1, \ldots , p \}$ belongs to 
$\ncdisc ( \, \{ 1, \ldots , p \} \, )$ and the partition
$\{ B_0'', B_1'', \ldots , B_t'' \}$ of 
$\{ p+1, \ldots , p+q \}$ belongs to 
$\ncdisc ( \, \{ p+1, \ldots , p+q \} \, )$.

$3^o$ If $\pi \in \ncannpq$, then all the $\card (B_0') \cdot \card (B_0'')$
permutations mentioned in part $1^o$ belong to $\sanc$; hence there are
exactly $\card (B_0') \cdot \card (B_0'')$ permutations $\tau \in \sanc$ 
such that $\orbits ( \tau ) = \pi$.

\vspace{10pt}

{\bf Proof.} $1^o$ This is because for every $(b',b'') \in B_0' \times B_0''$
there exists a unique $\tau \in \GSpq$ which is standard in $(p,q)$--annular 
sense, has orbits$( \tau ) = \pi$, and satisfies the relation 
$\tau ( b') = b''$. The formulas describing how this unique $\tau$ acts on
its orbits are:
$\tau \mid B_i ' = \gamma_{ext} \mid B_i '$, $1 \leq i \leq s$; 
$\tau \mid B_j '' = \gamma_{int} \mid B_j ''$, $1 \leq j \leq t$; 
and $\tau \mid B_0 \ = \ ( b_1', \ldots , b_k' , b_1'', \ldots , b_l'' ),$
where $b_1', \ldots , b_k'$ is the enumeration of $B_0'$ made such that
$\gamma_{ext} \mid B_0' = ( b_1', \ldots , b_k')$ and $b_k' = b'$, while
$b_1'', \ldots , b_l''$ is the enumeration of $B_0''$ made such that
$\gamma_{int} \mid B_0'' = ( b_1'', \ldots , b_l'')$ and $b_1'' = b''$.

$2^o$ ``$\Rightarrow$'' Suppose that $\pi \in \ncannpq$, and let 
$\tau \in \sanc$ be such that $\orbits ( \tau ) = \pi$. If the partition 
$\{ B_0', B_1', \ldots , B_s' \}$ (respectively
$\{ B_0'', B_1'', \ldots , B_t'' \}$) was not non-crossing in disc sense, 
then we would find 4 distinct elements $a,b,c,d \in \{ 1, \ldots , p \}$
(respectively $a,b,c,d \in \{ p+1, \ldots , p+q \}$) such that $\tau$
satisfies (AC-1) for $a,b,c,d$ -- contradiction.

``$\Leftarrow$'' Let $\tau$ be any of the $\card (B_0') \cdot \card (B_0'')$
permutations of $[p+q]$ which are standard in $(p,q)$--annular sense and 
have the blocks of $\pi$ as orbits. We will prove that $\tau \in \sanc$ 
(which will imply that $\pi = \orbits ( \tau ) \in \ncannpq$). The fact 
that $\tau$ does not satisfy (AC-1) is immediate from the hypothesis that 
the partitions $\{ B_0', B_1', \ldots , B_s' \}$ and 
$\{ B_0'', B_1'', \ldots , B_t'' \}$ are non-crossing in disc sense; so 
it remains to check that $\tau$ does not fulfill any of (AC-2), (AC-3).

Suppose that $a,b,c,x,y$ are distinct elements of $[p+q]$ for which $\tau$
satisfies (AC-2). The elements $x,y$ have to be from the two distinct orbits
of $\gamma$ -- say that $x \in \{ 1, \ldots , p \}$ and
$y \in \{ p+1, \ldots , p+q \}$. Since $x$ and $y$ must belong to the same 
orbit of $\tau$, and since the only $(p,q)$--connecting orbit of $\tau$ is 
$B_0$, it follows that $x \in B_0'$ and $y \in B_0''$. Then $a,b,c$ must all 
belong to one of the orbits $B_1', \ldots , B_s', B_1'', \ldots , B_t''$,
which means in particular that either $a,b,c \in \{ 1, \ldots , p \}$ or
$a,b,c \in \{ p+1, \ldots , p+q \}$. But this leads to a contradiction with 
the fact that $\tau$ is standard in $(p,q)$--annular sense, exactly as shown 
in (3.5) of Remark 3.6. So we conclude that $\tau$ does not satisfy (AC-2).

Suppose that $a,b,c,d,x,y$ are distinct elements of $[p+q]$ for which $\tau$
satisfies (AC-3). Exactly as in the preceding paragraph we see that (after 
swapping $x$ and $y$, if necessary) we have $x \in B_0'$ and $y \in B_0''$.
It cannot happen that $a,c \in B_i'$ and $b,d \in B_j''$ for some 
$1 \leq i \leq s$ and $1 \leq j \leq t,$ because it would follow that 
$a,c \in \{ 1, \ldots , p \}$, $b,d \in \{ p+1, \ldots , p+q \}$, and this is 
not compatible with the fact that 
$\lambda_{x,y} \mid \{ a,b,c,d \} = (a,b,c,d)$
(which is part of (AC-3)). The case when $a,c \in B_j''$ and $b,d \in B_i'$
(for some $1 \leq j \leq t$ and $1 \leq i \leq s$) cannot occur because of the
same reason. It remains that: either $a,c \in B_i'$, $b,d \in B_j'$ for some
$1 \leq i,j \leq s,$ $i \neq j$; or $a,c \in B_i''$, $b,d \in B_j''$ for some
$1 \leq i,j \leq t,$ $i \neq j$. But then we obtain that 
$\gamma_{ext} \mid \{ a,b,c,d \} = \lambda_{x,y} \mid \{ a,b,c,d \}$ = 
$(a,b,c,d)$ (respectively that 
$\gamma_{int} \mid \{ a,b,c,d \} = \lambda_{x,y} \mid \{ a,b,c,d \}$ = 
$(a,b,c,d)$), and we get a contradiction with the hypothesis that $B_i'$ 
and $B_j'$ (respectively $B_i''$ and $B_j''$) do not cross. So we conclude 
that $\tau$ does not satisfy (AC-3) -- hence that $\tau \in \sanc$.

$3^o$ This was proved at the same time with the part ``$\Leftarrow$'' of
$2^o$. {\bf QED}

$\ $

{\bf 4.7 Remark.} The relation between annular non-crossing partitions and
permutations isn't so good as the one we had in the disc case, because the 
natural map $\sanc \to \ncannpq$ is not one-to-one. Nevertheless, the failure
of the injectivity of this map is limited -- the only pathology that can 
appear is the one described in Proposition 4.6.

This pathology has a consequence which was observed in \cite{Ki}, concerning
the annular counterpart for the complementation map found by Kreweras on 
$\ncdisc (n)$. We will show that the annular version of the 
Kreweras complementation map is well-defined and bijective as a map from 
$\sanc$ to itself (see Remark 6.7 below). But when trying to define the 
Kreweras complementation as a map on $\ncannpq$, one runs into the problem 
(noticed in \cite{Ki}) that the map is not well-defined on the set of 
partitions $\pi \in \ncannpq$ which have exactly one $(p,q)$--connecting block,
and that it is not bijective on the part of $\ncannpq$ where it is defined.

$\ $

{\bf 4.8 Remark} {\em (annular non-crossing complete matchings).} 
We say that a partition $\pi$ of $[p+q]$ is a {\em complete matching}
if every block of $\pi$ has exactly 2 elements.

If $\pi$ is a complete matching of $[p+q],$ then there exists a unique 
permutation $\tau \in \GSpq$ such that $\orbits ( \tau ) = \pi$. This makes 
clear that the natural map from $\sanc$ to $\ncannpq$ becomes a bijection
when it is restricted to go from $\{ \tau \in \sanc \ | $ every orbit of 
$\tau$ has exactly $\mbox{2 elements} \}$ to $\{ \pi \in \ncannpq \ |$
$\pi$ is a $\mbox{complete matching} \}$.

If one only wants to work with annular non-crossing complete matchings, then 
it is easy to concoct a definition of these objects which is quite a bit 
simpler than what we had in Section 3. Essentially one stipulates that a
complete matching $\pi$ of $[p+q]$ is $(p,q)$--annular non-crossing if and only
if every group of up to 3 blocks ( = pairs) of $\pi$ can be drawn without 
crossings in the $(p,q)$--annulus. (Explanation: Let $\tau$ be the unique 
permutation of $[p+q]$ such that $\orbits ( \tau ) = \pi$. Then the conditions
(ANS-1), (ANS-2) and (AC-2) can't apply to $\tau$. The non-fulfillment of 
(AC-1) amounts to the fact that 2 blocks of $\pi$ which are both contained 
either in $\{ 1, \ldots , p \}$ or in $\{ p+1, \ldots , p+q \}$ do not cross.
The non-fulfillment of (AC-3) amounts to the fact that a group of 3 blocks of
$\pi$, out of which at least one is $(p,q)$--connecting, can always be drawn
without crossings in the $(p,q)$--annulus.)

Annular non-crossing complete matchings (and non-crossing complete matchings 
drawn in multi-annuli as well) have been used for a long time in the physics 
literature, under the name of planar Feynman diagrams. In the mathematical 
literature they can also be traced back quite a while, at least to the paper
\cite{Tu1} by Tutte. More recently, annular non-crossing complete matchings 
have been studied by Jones \cite{J1}, as part of a discussion on annular 
counterparts of the Temperley-Lieb algebras. 

$\ $

$\ $

\section{A relation between non-crossing permutations in disc sense and 
in annular sense}

In this section we continue to keep fixed the positive integers $p,q,$
\setcounter{equation}{0}
and the terminology related to the ``$(p,q)$--annulus'' which was 
introduced in the Section 3. In particular,
\[
\gamma_{ext},\ \gamma_{int}, \ \gamma, \ \lambda_{x,y} \in \GSpq
\]
will be exactly as in the Notations 3.2. On the other hand we will also use 
the notation
\begin{equation}
\gamma_o \ := \ (1,2, \ldots , p+q-1, p+q) \in \GSpq ,
\end{equation}
which comes from considerations on non-crossing permutations in 
the disc sense, as in Section 2 (the value of $n$ from Section 2 being 
now $n = p+q$).

At first glance one wouldn't be too tempted to relate the sets 
$\GS_{disc-nc} (p+q)$ (defined as in Section 2) and 
$\sanc$ (defined in Section 3), because one draws different types of
pictures for the permutations belonging to these two sets. Nevertheless,
it turns out that in our algebraic treatment we have the following simple 
relation between them:

$\ $

{\bf 5.1 Theorem.} $\sanc$ is the smallest subset of $\GSpq$ which contains 
\newline
$\GS_{disc-nc}(p+q)$ and which is invariant under conjugation with 
$\gamma_{ext}$ and with $\gamma_{int}.$

$\ $

The proof of the theorem relies on three facts which we state in a 
separate lemma.

$\ $

{\bf 5.2 Lemma.} 
$1^o$ $\sanc$ is invariant under conjugation with 
$\gamma_{ext}$ and with $\gamma_{int}$.

$2^o$ $\GS_{disc-nc} (p+q) \subset \sanc$.

$3^o$ Let $\tau \in \sanc$ be such that $\tau (p) = p+1.$
Then $\tau \in \GS_{disc-nc} (p+q).$

$\ $

{\bf Proof of Theorem 5.1} {\em (by using the Lemma 5.2).}
After taking into account the statements $1^o$ and $2^o$ of the lemma,
all that remains to be shown is that for every 
$\tau \in \sanc$ there exist $0 \leq u < p$, $0 \leq v <q$ such that:
\begin{equation}
\gamma_{ext}^u \, \gamma_{int}^v \, \tau \, \gamma_{int}^{-v} \,
\gamma_{ext}^{-u} \in \GS_{disc-nc} ( p+q).
\end{equation}
We pick a permutation $\tau \in \sanc$ about which we show that (5.2) holds.
If $\tau$ is $(p,q)$--disconnected, then the Remark 3.8 gives us that
$\tau \in \GS_{disc-nc} (p+q),$ and (5.2) holds with $u=v=0.$ So let 
us assume that $\tau$ is $(p,q)$--connected. Let $a \in \{ 1, \ldots , p \}$
and $b \in \{ p+1, \ldots , p+q \}$ be such that $\tau (a) = b.$ There exist
$0 \leq u < p$ and $0 \leq v <q$ such that $\gamma_{ext}^u (a) = p,$
$\gamma_{int}^v (b) = p+1;$ for these $u,v$ we set $\sigma :=
\gamma_{ext}^u \, \gamma_{int}^v \, \tau \, \gamma_{int}^{-v} \,
\gamma_{ext}^{-u}.$ We have that $\sigma \in \sanc$ by Lemma 5.2.1. On 
the other hand we have that $\sigma (p) = p+1$; so the Lemma 5.2.3
applies to give us that $\sigma \in \GS_{disc-nc} ( p+q)$, and (5.2) 
follows. {\bf QED}

$\ $

{\bf Proof of the Lemma 5.2.} We will show in detail the proof of the 
statement $3^o$. The proofs of $1^o$ and $2^o$ are similar in spirit
(quite straightforward, but a bit tedious), and we will leave them as 
an exercise to the meticulous reader.

For the less meticulous reader, we can point the following ``pictorial''
arguments supporting the statements $1^o$ and $2^o$. (These arguments
have some merit in view of the discussion in Section 3, which makes a 
case that the annular concepts of ``planar'' and ``non-crossing'' are
in fact identical.)

-- For the statement $1^o$: The operation of conjugating $\tau \in \sanc$ 
with a permutation of the form $\gamma_{ext}^u \gamma_{int}^v$ for some 
$0 \leq u <p$, $0 \leq v <q$ can be simply viewed as the operation of
changing cyclically the labels of the $p$ points on the external circle
and of the $q$ points on the internal circle of the $(p,q)$-annulus,
without touching the actual picture of $\tau$. Hence the resulting 
permutation $\tau ' = ( \gamma_{ext}^u \gamma_{int}^v ) \tau 
( \gamma_{ext}^u \gamma_{int}^v )^{-1}$ must also be in $\sanc$.

-- For the statement $2^o$: Let $\tau$ be a permutation in 
$\GS_{disc-nc} (p+q)$. Then $\tau$ has a non-crossing picture in the disc,
as discussed in the Remark 2.6.1. Instead of being drawn in a disc, 
this picture can (clearly) be also drawn in a square, such that the points 
$1, \ldots , p$ are marked on the top horizontal side of the square
(from left to right) and the points $p+1, \ldots , p+q$ are marked on 
the bottom horizontal side (from right to left). Then let us fold 
this square into a cylinder, by glueing together its vertical sides; and
after that let us flatten the resulting cylinder, turning it into the 
$(p,q)$-annulus. In the process, the disc non-crossing picture of $\tau$ 
is first turned into a picture drawn on the lateral part of the cylinder, 
and then becomes a $(p,q)$--annular picture (thus showing that 
$\tau \in \sanc$).

Let us mention once again that it it is quite easy (though perhaps tedious)
to also verify the statements $1^o$ and $2^o$ in a purely algebraic fashion, 
obtaining an argument similar in spirit to the one shown next for the 
statement $3^o$.

So, in order to verify the statement $3^o$ we fix (from now on and until the 
end of the proof) a permutation $\tau \in \sanc$ such that $\tau (p)= p+1$.
Our goal is to show that $\tau \in \GS_{disc-nc} (p+q)$.

Note that the permutation $\lambda_{p,p+1}$ (defined as in 
Notations 3.2) is $\lambda_{p,p+1}$ = 
$(1, \ldots p-1, p+2, \ldots , p+q )(p)(p+1) \in \GSpq$;
as a consequence, we have that:
\begin{equation}
\lambda_{p,p+1} \mid [p+q] \setminus \{p,p+1 \} \ = \ 
\gamma_o \mid [p+q] \setminus \{p,p+1 \} .
\end{equation}

A special role in the proof will be played by the orbit of $\tau$ which 
contains $p$ and $p+1$. We denote this orbit as $A$. Let us observe that 
the cycle $\tau \mid A$ is of the form
\begin{equation}
\left\{  \begin{array}{c}
\tau \mid A = ( a_1, \ldots , a_k, b_1, \ldots , b_l ), \ \mbox{ where} \\
1 \leq a_1 < \cdots < a_k = p < p+1 = b_1 < \cdots < b_l \leq p+q.
\end{array}  \right.
\end{equation}
Indeed, the Remark 3.4.1 (Case 2) shows that 
$\tau \mid A = ( a_1, \ldots , a_k, b_1, \ldots , b_l )$  where
$a_1, \ldots , a_k$ $\in \{ 1, \ldots , p \}$,
$b_1, \ldots , b_l \in \{ p+1, \ldots , p+q \}$, and where both 
$a_1, \ldots , a_k$ and $b_1, \ldots , b_l$ can be permuted cyclically 
to be put in increasing order. Since $a_k$ is the unique element of 
$A \cap \{ 1, \ldots , p \}$ mapped by $\tau$ into $\{ p+1, \ldots , p+q \}$,
we must have $a_k = p$ and $b_1 = \tau (a_k) = p+1.$ Moreover, if 
$a_k =p,$ then the only cyclic permutation of $a_1, \ldots , a_k$ which can 
possibly put these numbers in increasing order is the identity permutation; 
hence $a_1 < \cdots < a_k,$ and a similar argument shows that 
$b_1 < \cdots < b_l.$ 

Let us prove that $\tau$ is standard in disc sense, i.e. that 
it does not satisfy the condition (DNS) of Proposition 2.8. Suppose by
contradiction that we found 3 distinct elements $a,b,c$ belonging to an 
orbit $B$ of $\tau$ such that $\gamma_o \mid \{ a,b,c \} = (a,b,c)$ and
$\tau \mid \{ a,b,c \} = (a,c,b).$ Then $B$ is distinct from the special 
orbit $A$ of the preceding paragraph (since $\tau \mid A = \gamma_o \mid A$,
by (5.4)). But then 
$\{ a,b,c \} \subset [p+q] \setminus \{ p,p+1 \}$, hence
$\lambda_{p,p+1} \mid \{ a,b,c \} = \gamma_o \mid \{ a,b,c \} = (a,b,c)$
(by (5.3)). This shows that $\tau$ satisfies (AC-2) for the 5 elements
$a,b,c,p,p+1$, contradicting the hypothesis that $\tau \in \sanc$.

We next observe that if $B,C$ are two distinct orbits of $\tau$ such that 
$B \neq A, \ C \neq A,$ then $B$ and $C$ are not crossing in the disc sense.
Indeed, in the opposite case we could find distinct elements $b,b' \in B,$
$c,c' \in C,$ such that $\gamma_o \mid \{ b,b',c,c' \} = (b,c,b',c')$. Since
$\{ b,b',c,c' \} \subset [p+q] \setminus \{ p,p+1 \}$, the Equation (5.3)
would imply that $\lambda_{p,p+1} \mid \{ b,b',c,c' \} = (b,c,b',c')$,
and it would follow that $\tau$ satisfies (AC-3) for the 6 elements 
$b,b',c,c',p,p+1$ (contradiction).

So we are left to fix an orbit $B$ of $\tau$ such that $B \neq A$, and to 
prove that $A$ and $B$ are not crossing in the disc sense. Suppose this
is not true. Then we can find $b,b' \in B$ and $a \in A$ such that $b<a<b'$
and such that either $A \cap \{ 1, \ldots , b-1 \} \neq \emptyset$ or
$A \cap \{ b'+1, \ldots , p+q \} \neq \emptyset$.
We distinguish 4 possible cases:

\vspace{6pt}

\begin{tabular}{lcl}

{\em Case 1:} $b' < p.$  &  &   {\em Case 3:} $b<p, \ b' > p+1,$ and 
                           $A \cap \{ 1, \ldots , b-1 \} \neq \emptyset$. \\
                         &  &                                             \\
{\em Case 2:} $b > p+1.$ &  &  {\em Case 4:} $b<p, \ b' > p+1,$ and 
                           $A \cap \{ b'+1, \ldots , p+q \} \neq \emptyset$.
\end{tabular}

\vspace{6pt}
\noindent
But now, each of these 4 cases comes in contradiction with the hypothesis that 
$\tau \in \sanc$. Indeed: In Case 1 we find that $\tau$ satisfies (AC-1) for 
$b,a,b',p$ and in Case 2 we find that $\tau$ satisfies (AC-1) for 
$p+1, b,a,b'$. In the Case 3 we pick an element 
$a' \in A \cap \{ 1, \ldots , b-1 \}$, and we find that $\tau$ satisfies (AC-2)
for $a',p,p+1,b,b'$ (this is because 
$\lambda_{b,b'} \mid \{ a',p,p+1 \} = (p,a',p+1)$, while (5.4) implies that 
$\tau \mid \{ a',p,p+1 \} = (a',p,,p+1)$). The Case 4 is similar to Case 3:
we pick an element $a'' \in A \cap \{ b'+1, \ldots , p+q \}$, and we find 
that $\tau$ satisfies (AC-2) for $p,p+1,a'',b,b'$. {\bf QED}

$\ $

$\ $

\section{The annular version of the geodesic condition} 

In this section we continue to maintain the notations introduced throughout
\setcounter{equation}{0}
the Sections 3--5 (in particular, the special permutations
$\gamma , \gamma_o \in \GSpq$ are as defined in Equations (3.1) and (5.1),
respectively). We will prove the following theorem:

$\ $

{\bf 6.1 Theorem.} $1^o$ Let $\tau$ be a permutation of $[p+q]$ which is 
$(p,q)$--connected. Then 
$\# ( \tau ) \, + \, \# ( \tau^{-1} \gamma ) \leq  p+q,$ and $\tau$ 
belongs to $\sanc$ if and only if the above inequality is an equality.

$2^o$ Let $\tau$ be a permutation of $[p+q]$ which is $(p,q)$--disconnected. 
Then $\# ( \tau ) \, + \, \# ( \tau^{-1} \gamma ) \leq  p+q+2,$ and $\tau$ 
belongs to $\sanc$ if and only if the above inequality is an equality.

$\ $

{\bf 6.2 Remark.} Before starting on the proof of Theorem 6.1, let us take
a moment to place its statement in the framework of the geodesic condition
of Section 2.10. Referring to the notations in that section, where we set
$n=p+q,$ we see that the equation
$d( \, id, \tau \, ) + d( \, \tau , \gamma \, ) = d( \, id, \gamma \, )$ 
is tantamount to
$\# ( \tau ) \, + \, \# ( \tau^{-1} \gamma ) = p+q+2$ (by Equation (2.6), and 
the fact that $\gamma$ has 2 orbits). So in the part $2^o$ of Theorem 6.1
we obtain precisely that $\tau$ lies on a geodesic between $id$ and
$\gamma$. In part $1^o$ of the theorem we get that the triangle inequality
for $id, \, \tau$ and $\gamma$ is strict:
$d( \, id, \tau \, ) + d( \, \tau , \gamma \, ) = d( \, id, \gamma \, )+2.$ 
It is easy to see from parity considerations that the equality 
$d( \, id, \tau \, ) + d( \, \tau , \gamma \, ) = d( \, id, \gamma \, )+1$ 
can never occur, so a suggestive way of stating the part $1^o$ of Theorem 6.1
would be thus: ``A $(p,q)$--connected permutation of $[p+q]$ is in $\sanc$
if and only if {\em it barely fails} Biane's geodesic condition''.

$\ $

In the proof of Theorem 6.1 we will use the following well-known fact (which 
also lies at the basis of the proof of the geodesic condition in the disc 
case -- see \cite{B1}, Lemma 1): for every $1 \leq i<j \leq p+q$ and every 
$\tau \in \GSpq$ we have
\begin{equation}
\# ( \, \tau \cdot (i,j) \, ) \ = \ \left\{  \begin{array}{cl}
\# ( \tau ) + 1  &  \mbox{if $i,j$ are in the same orbit of $\tau$}  \\
\# ( \tau ) - 1  &  \mbox{otherwise.}
\end{array}  \right.
\end{equation}

We will also use the following observation:

$\ $

{\bf 6.3 Lemma.} Let $\G$ denote the subgroup of $\GSpq$ which is generated
by $\gamma_{ext}$ and $\gamma_{int}$. Let $\U$ be a set of $(p,q)$--connected
permutations of $[p+q]$, such that $\U$ is invariant under conjugation by 
elements of $\G$. Denote $\U_o := \{ \tau \in \U \ | \ \tau (p) = p+1 \}$. 
Then
$\U = \{ \sigma \tau \sigma^{-1} \ | \ \sigma \in  \G, \ \tau \in \U_o \}$.

$\ $

The easy proof of Lemma 6.3 is similar to the argument used in the proof of 
Theorem 5.1, and is left to the reader.

$\ $

{\bf 6.4 Lemma.} Let $\tau$ be a permutation of $[p+q]$ such that 
$\tau (p) = p+1.$ Then $\tau \in \sanc$ if and only if
$\# ( \tau ) \, + \, \# ( \tau^{-1} \gamma ) = p+q.$

\vspace{10pt}

{\bf Proof.} It is immediate that $\gamma$ and $\gamma_o$ are related by
\begin{equation}
\gamma_o \ = \ \gamma \cdot (p,p+q).
\end{equation}
The hypothesis that $\tau (p) = p+1$ implies that 
$( \tau^{-1} \gamma ) (p+q) = \tau^{-1} (p+1) =p;$ hence $p$ and $p+q$
belong to the same orbit of $\tau^{-1} \gamma$, and we get:
\begin{equation}
\# ( \tau^{-1} \gamma_0 )  \begin{array}[t]{lll}
= & \# ( \ \tau^{-1} \gamma \cdot (p,p+q) \ ) & \mbox{ (by (6.2))}   \\
= & \# ( \tau^{-1} \gamma ) +1 & \mbox{ (by (6.1)).} 
\end{array}
\end{equation}
But then, for the permutation $\tau$ given in the lemma we can write the 
following equivalences:
\[
\begin{array}{cl}
\tau \in \sanc               &                                      \\
\Updownarrow                 &                                      \\
\tau \in \GS_{disc-nc} (p+q) & \mbox{ (by the Lemma 5.2)}  \\
\Updownarrow                 &                                      \\
\# ( \tau ) \, +  \, \# ( \tau^{-1} \gamma_o ) = p+q+1
                  & \mbox{ (geodesic condition in the disc case)}   \\
\Updownarrow                 &                                      \\
\# ( \tau ) \, +  \, \# ( \tau^{-1} \gamma ) = p+q
                             & \mbox{ (by (6.3)). {\bf QED} }
\end{array}
\]

$\ $

{\bf Proof of Theorem 6.1.} 
$2^o$ As observed in the Remark 6.2, the equation 
$\# ( \tau ) \, + \, \# ( \tau^{-1} \gamma ) = p+q+2$ can be written as
$d( \, id, \tau \, ) + d( \, \tau , \gamma \, ) = d( \, id , \gamma \, ).$
By the results from the disc case (in the more general form reviewed 
in (2.8) of Section 2.10), the fulfillment of this equation is equivalent  
to the fact that $\{ 1, \ldots , p \}$ and $\{ p+1, \ldots , p+q \}$ are
$\tau$-invariant, and that 
$\tau \mid \{ 1, \ldots ,p \} \in \GS_{disc-nc}( \, \{ 1, \ldots , p \} \, )$,
$\tau \mid \{ p+1, \ldots ,p+q \} \in \GS_{disc-nc}( \, \{ p+1, \ldots , 
p+q \} \, )$. Finally, the latter fact is equivalent to the statement that
$\tau \in \sanc$, by (3.6) of Remark 3.8. 

The inequality
$\# ( \tau ) \, + \, \# ( \tau^{-1} \gamma )  \leq p+q+2$ holds for all 
$\tau \in \GSpq$ because it is a reformulation of the triangle inequality
$d( \, id, \tau \, ) + d( \, \tau , \gamma \, ) \geq d( \, id , \gamma \, ).$

$1^o$  If $\tau \in \GSpq$ is $(p,q)$--connected then the joint action of 
$\tau$ and $\gamma$ on $[p+q]$ is transitive, and the equality (2.10) from 
Section 2.11 (used for $\tau$ and $\gamma$, with $\# ( \gamma ) =2$ and 
$\# ( \tau \vee \gamma ) = 1$) gives us that
$\# ( \tau ) \, + \, \# ( \tau^{-1} \gamma ) \leq  p+q$.

Let us denote
$\U ' := \{ \tau \in \sanc \ | \ \tau$ is $(p,q)$--connected$\}$ and
$\U '' := \{ \tau \in \GSpq \ | \ \tau$ is $(p,q)$--connected and
$\# ( \tau ) \, + \, \# ( \tau^{-1} \gamma ) = p+q \} .$
It is immediately verified that both $\U'$ and $\U''$ satisfy the hypotheses
of Lemma 6.3. Hence the Lemma 6.3 will give us that $\U ' = \U '',$ if we can
prove that $\U_o ' = \U_o ''$, where 
$\U_o ' :=  \{ \tau \in \U ' \ | \ \tau (p) = p+1 \}$ and
$\U_o '' :=  \{ \tau \in \U '' \ | \ \tau (p) = p+1 \}$.
But the equality of $\U_o '$ and $\U_o ''$ is precisely the statement of 
Lemma 6.4.  {\bf QED}

$\ $

{\bf 6.6 Corollary.} If $\tau \in \sanc$, then $\tau^{-1} \gamma \in \sanc$.

\vspace{10pt}

{\bf Proof.} Note that $\tau^{-1} \gamma$ is $(p,q)$--connected if and only if 
$\tau$ is so. (This is simply because a product of two $(p,q)$--disconnected 
permutations is again $(p,q)$--disconnected.) Since, as is easily seen, the 
quantity $\# ( \tau ) \, + \, \# ( \tau^{-1} \gamma )$ does not change when 
we replace $\tau$ by $\tau^{-1} \gamma$, the statement of the corollary 
follows from the one of Theorem 6.1. {\bf QED}

$\ $

{\bf 6.7 Remark.} Due to the above corollary, it makes sense to define a map
$K$ from $\sanc$ to itself by setting $K( \tau ) := \tau^{-1} \gamma$,
$\tau \in \sanc$. It is clear that $\tau$ is injective, so (being a self-map
of a finite set) it has to be a bijection. By analogy with the disc case, 
we will term $K$ as {\em ``the annular Kreweras complementation map''.} The 
name is justified by the fact that $K$ has a pictorial description which 
parallels the original construction made by Kreweras \cite{K} in the disc case.
We briefly describe how this goes, in the more interesting situation when 
$\tau$ is $(p,q)$--connected (if $\tau \in \sanc$ is $(p,q)$--disconnected, 
then $K( \tau )$ is obtained by taking separately the Kreweras complements 
in disc sense for $\tau \mid \{ 1, \ldots , p \}$ and for 
$\tau \mid \{ p+1, \ldots , p+q \}$). We proceed as follows: On the external 
circle of the $(p,q)$--annulus we mark $p$ new points labelled 
$\overline{1}, \ldots , \overline{p}$, such that $\overline{1}$ lies between 
1 and 2, $\overline{2}$ lies between $2 \mbox{ and } 3, \ldots , \overline{p}$ 
lies between $p$ and 1. Similarly, on the internal circle we mark $q$ new
points labelled $\overline{p+1}, \ldots , \overline{p+q},$ such that 
$\overline{p+1}$ lies between $p+1$ and $p+2, \ldots , \overline{p+q}$ lies
between $p+q$ and $p+1.$ Let $\tau$ be a $(p,q)$--connected permutation in 
$\sanc$, and let us draw a $(p,q)$--annular planar picture of $\tau$ by using 
the points $1, \ldots , p+q$. Then  $K( \tau )$ is in some sense the 
``maximal'' permutation in $\sanc$ which can be drawn by using the points
$\overline{1}, \ldots , \overline{p+q},$ such that the admissible pictures
of the $p+q$ cycles of $\tau$ and $K( \tau )$ (taken together!) enclose 
regions which are pairwise disjoint.
See the Figure 10 at the end of the paper for a concrete example.

$\ $

Another corollary of Theorem 6.1 refers to the enumeration of $\sanc$. The
number of $(p,q)$--disconnected permutations in $\sanc$ is equal to 
$\Bigl( (2p)! (2q)! \Bigr)/ \Bigl( p! (p+1)! q! (q+1)! \Bigr)$. This is 
immediate from the Remark 3.8 and the well-known fact that $\sdnc$ (or
equivalently $\ncdisc (n)$) is counted by the Catalan number
$(2n)! / n! (n+1)!,$ $n \geq 1.$ In the $(p,q)$--connected case we have:

$\ $

{\bf 6.8 Corollary.} The number of $(p,q)$--connected permutations in $\sanc$
is:
\begin{equation}
\frac{2pq}{p+q} \cdot 
\left( \begin{array}{c} 2p-1 \\ p \end{array} \right) \cdot
\left( \begin{array}{c} 2q-1 \\ q \end{array} \right) .
\end{equation}

\vspace{10pt}

{\bf Proof.} The $(p,q)$--connectedness of a permutation $\tau \in \GSpq$ 
is equivalent to the fact that the joint action of $\tau$ and $\gamma$ is 
transitive on $[p+q]$. Re-denoting $\tau =: \sigma_1$, 
$\tau^{-1} \gamma =: \sigma_2$ we thus find, in view
of Theorem 6.1, that the number of $(p,q)$--connected permutations in 
$\sanc$ is equal to the number of couples $( \sigma_1 , \sigma_2 )$ such that:
$\sigma_1 , \sigma_2 \in \GSpq$, $\sigma_1  \sigma_2 = \gamma$,
$\# ( \sigma_1 ) + \# ( \sigma_2 ) = p+q$, and the group generated by 
$\sigma_1$ and $\sigma_2$ acts transitively on $[p+q]$. The latter number is 
known to have the expression stated in (6.4) -- see the paper \cite{BMS} by 
Bousquet-M\'elou and Schaeffer, Theorem 1.2. {\bf QED}

$\ $

{\bf 6.9 Remarks.} $1^o$ The result from \cite{BMS} cited above is obtained 
via the enumeration of some planar maps called ``constellations''; so in 
a certain sense, the discussion has returned here to a class of planar maps
(even though these are not the same as the ones introduced back in 
Section 3.1, under the name of ``pictures of permutations'').

$2^o$ Let $\sancpair$ denote the set of permutations $\sigma \in \stanc$
such that every orbit of $\sigma$ has exactly 2 elements. It has been known
since long ago (see formula (1.1) in \cite{Tu1}) that the number of 
$(2p,2q)$--connected permutations in $\sancpair$ is precisely the double of
the number appearing in (6.4).

Let us say that a permutation $\sigma$ of $[2p+2q]$ is parity-alternating
if $\sigma (i) -i$ is odd for every $1 \leq i \leq 2(p+q)$. Let $\A$ denote
the set of $(2p,2q)$--connected permutations in $\sancpair$ which are 
parity-alternating. It is easily seen that $\A$ contains exactly one half of
the $(2p,2q)$--connected permutations belonging to $\sancpair$. Hence the 
coincidence observed in the preceding paragraph amounts to the fact that
$\A$ has the same number of elements as the set $\B$ of $(p,q)$--connected 
permutations in $\sanc$. We leave it as an amusing exercise to the reader to 
verify that a natural bijection between $\A$ and $\B$ can be defined in terms 
of annular Kreweras complementation maps, as introduced in Remark 6.7. The 
formula for this bijection (going from $\B$ to $\A$) is
\begin{equation}
\B \ni \tau \mapsto \widetilde{K} \Bigl( \,
\tau^{(odd)} \cup K( \tau )^{(even)} \, \Bigr) \in \A , 
\end{equation}
where the following notations are used:
\newline
$\bullet$ $K$ is the Kreweras complementation map on $\sanc$ (same as in 
Remark 6.7).
\newline
$\bullet$ $\widetilde{K}$ is the Kreweras complementation map on 
$\stanc$ (this is the map $\sigma \mapsto \sigma^{-1} \widetilde{\gamma},$
where $\widetilde{\gamma} := (1, \ldots , 2p)(2p+1, \ldots , 2p+2q)$).
\newline
$\bullet$ If $\tau$ is a permutation of $[p+q],$ then $\tau^{(odd)}$
and $\tau^{(even)}$ are the permutations of $\{ 1,3, \ldots ,$
$2(p+q)-1 \}$ and respectively of $\{ 2,4, \ldots , 2(p+q) \}$ defined by
$\tau^{(odd)} (2i-1) = 2 \tau (i) -1$ and 
$\tau^{(even)} (2i) = 2 \tau (i)$, $1 \leq i \leq p+q$.
\newline
$\bullet$ If $\tau_1$ and $\tau_2$ are permutations of 
$\{ 1,3, \ldots , 2(p+q)-1 \}$ and respectively $\{ 2,4, \ldots , 2(p+q) \}$,
then $\tau_1 \cup \tau_2$ denotes the permutation of $[ 2p+2q ]$ obtained
by joining $\tau_1$ and $\tau_2$ together.

$\ $

$\ $

\section{Relation with random matrices}

In this section we show how annular non-crossing permutations appear 
\setcounter{equation}{0}
in connection to second order asymptotics for certain random matrices. 
We will illustrate the phenomenon on a family of complex Wishart matrices 
(cf. Sections 7.3, 7.4 below). At the end of the section 
we will comment on how the same phenomenon also appears in connection to
Gaussian Hermitian matrices; but this is not so illustrative for 
our purposes here, as that example only involves complete matchings, rather 
than dealing with general permutations.

$\ $

{\bf 7.1 Framework.} $1^o$ Throughout this section $( \Omega , \F , P)$ is 
a fixed probability space, over which our random variables (measurable 
functions $f: \Omega \to \C$) will be considered. We will only deal with 
random variables which have finite moments of all orders, i.e. belong to:
\[
\LL_{-}^{\infty} ( \Omega , \F , P) \ := \ \Bigl\{ f: \Omega \to \C 
\begin{array}{lll}
\vline & f \mbox{ measurable }    &                                \\
\vline & \int | f( \omega ) |^n \ d P( \omega ) < \infty , & 
         \forall \ n \geq 1 
\end{array}  \Bigr\} .
\]
$\LL_{-}^{\infty} ( \Omega , \F , P)$ is a unital algebra of functions,
which is closed under conjugation, and is endowed with an ``expectation''
linear functional $\E$ defined by
\[
\E (f) \ := \  \int f( \omega ) \ d P( \omega ), \ \ 
f \in \LL_{-}^{\infty} ( \Omega , \F , P).
\]

$2^o$ Matrices over the algebra $\LL_{-}^{\infty} ( \Omega , \F , P)$
will be called random matrices over $( \Omega , \F , P)$. For a square 
random matrix $A = [ f_{i,j} ]_{i,j=1}^N$ over  $( \Omega , \F , P)$, the
normalized trace of $A$ is 
\[
\tr (A) \ := \ \frac{1}{N} \sum_{i=1}^N f_{i,i} \in
\LL_{-}^{\infty} ( \Omega , \F , P).
\]

$3^o$ We will work with independent families of standard complex 
Gaussian random variables (in the sense used e.g. in \cite{J}, page 13,
and described explicitly in the next lemma).
Our discussion will only involve the combinatorics of the moments 
of such a family; more precisely, we will only use the facts that the 
family is in $\LL_{-}^{\infty} ( \Omega , \F , P)$, and that it obeys 
the formulas (7.1), (7.2) given below. (These formulas are well-known, 
they are a version of the so-called ``Wick's Lemma'' -- cf. the
Sections 1.3 and 1.4 of the monograph \cite{J}.)

$\ $

{\bf 7.2 Lemma} {\em (Wick).} Let 
$\{ f_{\lambda} \ | \ \lambda \in \Lambda \}$ be complex
random variables over $( \Omega , \F , P)$ such that the real and the 
imaginary part of every $f_{\lambda}$ are centered Gaussian variables of
variance $1/2$, and such that the family 
$\{ \mbox{Re}(f_{\lambda}) | \lambda \in \Lambda \} \cup
\{ \mbox{Im}(f_{\lambda}) | \lambda \in \Lambda \}$ is independent.

$1^o$ Let $m,n$ be positive integers, $m \neq n,$ and consider two functions 
$\alpha : [m] \to \Lambda$ and $\beta : [n] \to \Lambda$. Then
\begin{equation}
\E \Bigl( \, f_{\alpha (1)} \cdots f_{\alpha (m)} \overline{f}_{\beta (1)}
\cdots \overline{f}_{\beta (n)} \, \Bigr) \ = \ 0.
\end{equation}

$2^o$ Let $n$ be a positive integer and consider two functions 
$\alpha, \beta : [n] \to \Lambda$. Then
\begin{equation}
\E \Bigl( \, f_{\alpha (1)} \cdots f_{\alpha (n)} \overline{f}_{\beta (1)}
\cdots \overline{f}_{\beta (n)} \, \Bigr) \ =  \
\mbox{card} \Bigl\{  \tau \in \GSn \ | \ \alpha = \beta \circ \tau \Bigr\} .
\end{equation}

$\ $

We can now introduce the special random matrices that we want to work with.

$\ $

{\bf 7.3 Notations.} In what follows, $G_1, \ldots , G_s$ ($s \geq 1$) 
will be a family of random $M \times N$ matrices over $( \Omega , \F , P)$,
with independent complex $N(0,1)$ entries. That is, we have:
\[
G_r \ = \ [ f_{i,j;r} ]_{ 1 \leq i \leq M, \ 1 \leq j \leq N} , 
\ \ \mbox{ for } 1 \leq r \leq s,
\]
where $\{ f_{i,j;r} \ | \ 1 \leq i \leq M, \ 1 \leq j \leq N, \
1 \leq r \leq s \}$ is an independent family of complex $N(0,1)$ random
variables. We denote
\[
X_r \ := \ \frac{1}{N} G_r^* G_r, \ \ \mbox{ for } 1 \leq r \leq s.
\]
If $w$ is an $n$-letter word over the alphabet $[s]$ (i.e. it is a function 
$w: [n] \to [s]$) for some $n \geq 1,$ then we will denote
\[
X_w \ := \ X_{w(1)} \cdots X_{w(n)}.
\]

When discussing  asymptotics, we will let $M$ and $N$ become the general 
terms of two sequences $(M_k)_{k=1}^{\infty}$ and $(N_k)_{k=1}^{\infty}$ which
increase to infinity in such a way that the limit 
$c = \lim_{k \to \infty} M_k/N_k$ exists and is in $(0, \infty ).$
For $M=M_k$ and $N=N_k$, the matrices $X_1, \ldots ,X_s$ and a generic
monomial $X_w$ made with them will be re-denoted as 
$X_1^{(k)}, \ldots , X_s^{(k)}$ and respectively as $X_w^{(k)}$.

$\ $

{\bf 7.4 Remark.} $X_1, \ldots , X_s$ are a particular case of
{\em complex Wishart} matrices. It is well-known that, when $k \to \infty$,
each of $X_1^{(k)}, \ldots , X_s^{(k)}$ converges in distribution to the 
Marcenko-Pastur distribution (for a nice presentation of this, see the 
Section 6 of \cite{HT}). The asymptotic joint behaviour of 
$X_1^{(k)}, \ldots , X_s^{(k)}$ is also remarkable: in the language of free
probability (see e.g. \cite{VDN}) we say that $X_1^{(k)}, \ldots , X_s^{(k)}$ 
are {\em asymptotically free} for $k \to \infty$. (For a recent discussion of 
this, see \cite{CC}.) Without entering into any details, we mention that both 
the asymptotic freeness of $X_1^{(k)}, \ldots , X_s^{(k)}$ and their individual
convergence to the Marcenko-Pastur distribution can be captured in a single 
formula, if one uses the theory of non-crossing cumulants of Speicher (see 
e.g. \cite{S}). The formula is:
\begin{equation}
\lim_{k \to \infty} \E \Bigl(  \tr  ( X_{w}^{(k)} ) \Bigr)  \ = \ 
\sum_{ \begin{array}{c}
{\scriptstyle \pi \in \ncdisc (n) \ such} \\
{\scriptstyle that \ w \ is \ constant} \\
{\scriptstyle on \ every \ block \ of \ \pi}
\end{array} } \ c^{ \ no. \ of \ blocks \ of \ \pi} ,
\end{equation}
where $n$ is a positive integer and $w$ is a word of length $n$ over the 
alphabet $[s]$.

The new fact that we want to put into evidence is that $\sanc$ shows 
up when the second-order asymptotics is considered. The quantities to look
at are of the form $\E \Bigl( \tr (X_v) \cdot \tr (X_w) \Bigr)$, where 
$v$ and $w$ are words over the alphabet $[s]$. (In the case $s=1$, when 
we deal with only one Wishart matrix, these quantities are related to the 
moments of the so-called ``2-point correlation function'' for the 
eigenvalues of the matrix.) The result is the following:

$\ $

{\bf 7.5 Theorem.} Let $v$ and $w$ be words of length $p$ and respectively
$q$ over the alphabet $[s]$. Then
\begin{equation}
\lim_{k \to \infty} 
\E \Bigl( \tr ( X_v^{(k)} ) \cdot \tr ( X_w^{(k)} ) \Bigr) \ = \
\lim_{k \to \infty} \E \Bigl( \tr ( X_v^{(k)} ) \Bigr) \cdot
\lim_{k \to \infty} \E \Bigl( \tr ( X_w^{(k)} ) \Bigr) 
\end{equation}
(where the two limits on the right-hand side of (7.4) can be described as
in (7.3)). Moreover, the sequence
$\E \Bigl( \tr ( X_v^{(k)} ) \cdot \tr ( X_w^{(k)} ) \Bigr) -
\E \Bigl( \tr ( X_v^{(k)} ) \Bigr) \cdot
\E \Bigl( \tr ( X_w^{(k)} ) \Bigr)$ 
goes to zero with an order of magnitude of $1/N_k^2$, and has
\[
\lim_{k \to \infty} N_k^2 \Bigl( \
\E ( \, \tr ( X_v^{(k)} ) \cdot \tr ( X_w^{(k)} ) \, ) -
\E ( \, \tr ( X_v^{(k)} ) \, ) \cdot
\E ( \, \tr ( X_w^{(k)} ) \, ) \  \Bigr)
\]
\begin{equation}
= \ \sum_{ \begin{array}{c}
{\scriptstyle \tau \in \sanc \ such \ that} \\
{\scriptstyle \tau \ is \ (p,q)-connected} \\
{\scriptstyle and \ (v \cup w) \circ \tau = v \cup w}
\end{array} } \ \ c^{ \# ( \tau ) } ,
\end{equation}
where $v \cup w$ denotes the juxtaposition of the words $v$ and $w$ (this
is a function from $[p+q]$ to $[s]$).

$\ $

The proof of Theorem 7.5 is based on the fact we have explicit summation 
formulas for the expectations appearing there, even before we let 
$k \to \infty$. This is explained in the following lemma (where we use the
notations $M,N, X_w$, without the extra index $k$).

$\ $

{\bf 7.6 Lemma.} $1^o$ Let $w$ be a word of length $n$ over the 
alphabet $[s]$. Then
\begin{equation}
\E \Bigl( \tr (X_w) \Bigr) \ = \ 
\frac{1}{N^{n+1}} \sum_{ \begin{array}{c}
{\scriptstyle \tau \in \GSn \ such}  \\
{\scriptstyle that \ w \circ \tau = w}
\end{array} } \ M^{ \# ( \tau ) } N^{ \# ( \tau^{-1} \gamma_o ) },
\end{equation}
where the notations for permutations are as in the preceding sections,
and in particular $\gamma_o$ stands for the forward cycle 
$(1, \ldots , n-1,n) \in \GSn$.

$2^o$ Let $v,w$ be words of length $p$ and respectively $q$ over the 
alphabet $[s]$. Then
\begin{equation}
\E \Bigl( \tr (X_v) \cdot \tr (X_w) \Bigr) \ = \ 
\frac{1}{N^{p+q+2}} \sum_{ \begin{array}{c}
{\scriptstyle \tau \in \GSpq \ such}  \\
{\scriptstyle that \ (v \cup w) \circ \tau = v \cup w}
\end{array} } \ M^{ \# ( \tau ) } N^{ \# ( \tau^{-1} \gamma ) } ,
\end{equation}
where, same as in the Sections 3-6,
$\gamma := (1, \ldots , p)(p+1, \ldots , p+q) \in \GSpq$.

$\ $

The proof of Lemma 7.6 is a straightforward computation, based on 
Equation (7.2) of Lemma 7.2; we present it for the reader's convenience. 
Substantial generalizations of the result of the lemma are known, but the
framework commonly considered is the one when $s=1$ -- see Corollary 2.4 in 
\cite{HSS}, or Theorem 2 in \cite {GLM}.

$\ $

{\bf Proof of Lemma 7.6.} We will show the argument for part $2^o$, the one 
for $1^o$ is similar. We have that
\[
\tr (X_v) \ = \ \tr ( \, X_{v(1)} \cdots X_{v(p)} \, )
\ = \ \frac{1}{N^p} \tr ( \, G_{v(1)}^* G_{v(1)} \cdots 
G_{v(p)}^* G_{v(p)} \, )
\]
\begin{equation}
= \ \frac{1}{N^{p+1}} \sum_{ \begin{array}{cc}
{\scriptstyle 1 \leq i_1, \ldots , i_p \leq M} \\
{\scriptstyle 1 \leq j_1, \ldots , j_p \leq N} 
\end{array} } \ 
\overline{f}_{i_1,j_1;v(1)} f_{i_1,j_2;v(1)}
\overline{f}_{i_2,j_2;v(2)} f_{i_2,j_3;v(2)} \cdots
\overline{f}_{i_p,j_p;v(p)} f_{i_p,j_1;v(p)} .
\end{equation}
$\tr (X_w)$ has a similar explicit formula, which we find convenient to 
write as:
\begin{equation}
\frac{1}{N^{q+1}} \sum_{ \begin{array}{cc}
{\scriptstyle 1 \leq i_{p+1}, \ldots , i_{p+q} \leq M} \\
{\scriptstyle 1 \leq j_{p+1}, \ldots , j_{p+q} \leq N} 
\end{array} } \ 
\overline{f}_{i_{p+1},j_{p+1};w(1)} f_{i_{p+1},j_{p+2};w(1)} \cdots
\overline{f}_{i_{p+q},j_{p+q};w(q)} f_{i_{p+q},j_{p+1};w(q)} .
\end{equation}
By multiplying together (7.8) and (7.9) we get a formula for 
$\tr (X_v) \cdot \tr (X_w)$, which is written more concisely if 
we record the indices $i_1, \ldots , i_{p+q}$ as a function 
$I:[p+q] \to [M]$ and the indices $j_1, \ldots , j_{p+q}$ as a function 
$J:[p+q] \to [N]$. It is also convenient to
use the word $v \cup w$, and replace $w(1), \ldots , w(q)$ by 
$(v \cup w) (p+1), \ldots , (v \cup w) (p+q).$ We obtain:
\begin{equation}
\tr (X_v) \cdot \tr (X_w) \ =
\end{equation}
\[
\frac{1}{N^{p+q+2}} \sum_{ \begin{array}{c}
{\scriptstyle I: [p+q] \to [M]} \\
{\scriptstyle J: [p+q] \to [N]} 
\end{array} } \ \Bigl( \prod_{m=1}^{p+q} 
\overline{f}_{I(m),J(m); (v \cup w)(m)} \Bigr) \cdot
\Bigl( \prod_{m=1}^{p+q} 
f_{I(m),J( \gamma (m)); (v \cup w)(m)} \Bigr) .
\]
We next apply $\E$ to both sides of (7.10), and use the Lemma 7.2 (in the
context where $\Lambda = [M] \times [N] \times [s],$ and where
$\alpha, \beta : [p+q] \to \Lambda$ are 
$\alpha (m) = (I(m),J( \gamma (m)), (v \cup w)(m))$,
$\beta (m) = (I(m),J(m), (v \cup w)(m))$). This gives us:
\begin{equation}
\E \Bigl(  \tr (X_v) \cdot \tr (X_w) \Bigr) \ = \ 
\end{equation}
\[
\frac{1}{N^{p+q+2}} \sum_{ \begin{array}{c}
{\scriptstyle I: [p+q] \to [M]} \\
{\scriptstyle J: [p+q] \to [N]} 
\end{array} } \ \card \{ \tau \in \GSpq \ | \ I \circ \tau = I, \
J \circ \tau = J \circ \gamma, \ (v \cup w) \circ \tau = v \cup w \} .
\]
The sum on the right-hand side of (7.11) can be re-written as a summation 
over $\tau$, namely:
\begin{equation}
\sum_{ \begin{array}{c}
{\scriptstyle \tau \in \GSpq \ such} \\
{\scriptstyle that \ (v \cup w) \circ \tau = v \cup w}
\end{array} } \ \card \Bigl\{  (I,J) \ \begin{array}{lll}
\vline & I:[p+q] \to [M],  &  I \circ \tau = I   \\
\vline & J:[p+q] \to [N],  &  J \circ \tau = J \circ \gamma   
\end{array} \Bigr\} .
\end{equation}
In (7.12), the condition $I \circ \tau = I$ is equivalent to asking that
$I$ is constant on the orbits of $\tau$; for any given $\tau$, there are 
$M^{ \# ( \tau )}$ ways of choosing such a function $I:[p+q] \to [M].$
Similarly, the condition $J \circ \tau = J \circ \gamma$ is equivalent to
$J \circ ( \gamma \tau^{-1} ) = J,$ hence to the fact that 
$J$ is constant on the orbits of $\gamma \tau^{-1}$. For any given 
$\tau \in \GSpq$, there are 
$N^{ \# ( \gamma \tau^{-1} )} = N^{ \# ( \tau^{-1} \gamma )}$ 
ways of choosing such a function $J:[p+q] \to [N].$ Thus we see that the
sum in (7.12) equals
$\sum_{ \begin{array}{c}
{\scriptstyle \tau \in \GSpq \ such} \\
{\scriptstyle that \ (v \cup w) \circ \tau = v \cup w}
\end{array} } \ M^{ \# ( \tau )} N^{ \# ( \tau^{-1} \gamma ) },$
and the formula (7.7) is obtained. {\bf QED}

$\ $

{\bf Proof of Theorem 7.5.} The Lemma 7.6.1 gives us that, for every 
$k \geq 1$:
\begin{equation}
\E \Bigl( \tr (X_v^{(k)} ) \Bigr) \ = \ \frac{1}{N_k^{p+1}} 
\sum_{ \begin{array}{c}
{\scriptstyle \tau_1 \in \GS  ( \, [p] \, ) \ such} \\
{\scriptstyle that \ v \circ \tau_1 = v}
\end{array} } \ M_k^{ \# ( \tau_1 ) } \cdot
N_k^{ \# ( \tau_1^{-1} \gamma_{ext} ) } ,
\end{equation}
\begin{equation}
\E \Bigl( \tr (X_w^{(k)} ) \Bigr) \ = \ \frac{1}{N_k^{q+1}} 
\sum_{ \begin{array}{c}
{\scriptstyle \tau_2 \in \GS  ( \, \{ p+1, \ldots , p+q \} \, ) } \\
{\scriptstyle such \ that \ w \circ \tau_2 ' = w}
\end{array} } \ M_k^{ \# ( \tau_2 ) } \cdot
N_k^{ \# ( \tau_2^{-1} \gamma_{int} )  } ,
\end{equation}
where in (7.14) we denoted by $\tau_2 '$ the permutation of $[q]$
corresponding to $\tau_2$ ($\tau_2 '$ sends $m$ to $\tau_2 (p+m) -p,$
for $1 \leq m \leq q$). We multiply together the Equations (7.13) and 
(7.14), and in the right-hand side of their product we put together 
$\tau_1$ and $\tau_2$ to form a permutation $\tau$ of $[p+q]$. We obtain:
\begin{equation}
\E \Bigl( \tr (X_v^{(k)}) \Bigr) \cdot 
\E \Bigr( \tr (X_w^{(k)}) \Bigr)  \ = \
\frac{1}{N_k^{p+q+2}} \sum_{ \begin{array}{c}
{\scriptstyle \tau \in \GSpq \ such \ that} \\
{\scriptstyle \tau \ is \ not \ (p,q)-connected} \\
{\scriptstyle and \ (v \cup w) \circ \tau = v \cup w}
\end{array} } \ \ M_k^{ \# ( \tau ) } 
N_k^{ \# ( \tau^{-1} \gamma ) } .
\end{equation}
When we subtract this out of the formula for 
$\E \Bigl( \tr ( X_v^{(k)} ) \cdot \tr ( X_w^{(k)} ) \Bigr)$ given by (7.7),
we get:
\[
N_k^2 \cdot \Bigl( \ \E ( \, \tr ( X_v^{(k)} ) \cdot \tr (  X_w^{(k)} ) \, ) 
- \E ( \, \tr (X_v^{(k)}) \, ) \cdot \E ( \, \tr (X_w^{(k)}) \, ) \ \Bigr)
\]
\begin{equation}
= \ \sum_{ \begin{array}{c}
{\scriptstyle \tau \in \GSpq \ such \ that} \\
{\scriptstyle \tau \ is \ (p,q)-connected} \\
{\scriptstyle and \ (v \cup w) \circ \tau = v \cup w}
\end{array} } \ \ (M_k/N_k)^{ \# ( \tau ) } \cdot
N_k^{ \# ( \tau ) + \# ( \tau^{-1} \gamma ) - (p+q) }.
\end{equation}
Finally, the Theorem 6.1 tells us that a $(p,q)$--connected permutation 
$\tau$ has $\# ( \tau ) + \# ( \tau^{-1} \gamma ) \leq p+q$, with equality 
if and only if $\tau$ is in $\sanc$; hence making $k \to \infty$ in (7.16)
leads us to (7.5). {\bf QED}

$\ $

{\bf 7.7 Remark.} Let us briefly point out that facts similar to those 
described in 7.4 -- 7.6 also hold when instead of Wishart matrices one uses
a family of Gaussian Hermitian (also called ``GUE'' -- cf. \cite{M}) random 
matrices with independent entries. In the setting of the Notations 7.3, such
a family $Y_1, \ldots , Y_s$ is obtained if we assume that $M=N$ (so that 
$G_1, \ldots , G_s$ are square matrices), and we define
\[
Y_r \ := \ \frac{1}{\sqrt{2N}} (G_r + G_r^* ), \ \ 1 \leq r \leq s.
\]
If $w$ is a word of length $n$ over the alphabet $[s]$, then we will denote
\[
Y_w \ := \ Y_{w(1)} Y_{w(2)} \cdots Y_{w(n)}.
\]

The counterpart of Equation (7.6) from Lemma 7.6.1 is the following formula
(see e.g. Theorem 3.1 in \cite{Th}): 
\begin{equation}
\E \Bigl( \tr (Y_w) \Bigr) \ = \ \frac{1}{N^{n+1}} \sum_{ 
\begin{array}{c}
{\scriptstyle \tau \ complete \ matching \ of \ [n]} \\
{\scriptstyle such \ that \ w \circ \tau = w}
\end{array} } \ \ N^{ \# ( \tau ) + \# ( \tau^{-1} \gamma_o )} ,
\end{equation}
where $w$ is a word of length $n$ over the alphabet $[s]$, and where by a
complete matching of $[n]$ we understand a permutation $\tau \in \GSn$ such
that every orbit of $\tau$ has exactly 2 elements. (Besides that, the 
notations used in (7.17) are identical to those from (7.6).)

The Equation (7.17) can be generalized without difficulty to deal with 
a product of two traces. We get that
\begin{equation}
\E \Bigl( \tr (Y_v) \cdot \tr (Y_w) \Bigr) \ = \ \frac{1}{N^{p+q+2}} \sum_{ 
\begin{array}{c}
{\scriptstyle \tau \ complete \ matching \ of \ [p+q]} \\
{\scriptstyle such \ that \ (v \cup w) \circ \tau = v \cup w}
\end{array} } \ \ N^{ \# ( \tau ) + \# ( \tau^{-1} \gamma )} ,
\end{equation}
where now $v$ and $w$ are words of length $p$ and respectively $q$ over the 
alphabet $[s]$ (and where the notations are analogous to those in 
Equation (7.7) of Lemma 7.6.2).

Let us now make $N$ become the general term of a sequence 
$(N_k)_{k=1}^{\infty}$ which goes to infinity. When $N=N_k$, we re-denote 
the matrices $Y_1, \ldots , Y_s$ and a generic monomial $Y_w$ made with
them as $Y_1^{(k)}, \ldots Y_s^{(k)}$ and as $Y_w^{(k)}$, respectively.
Starting from (7.17) and (7.18), it is easy to derive the counterpart of
Theorem 7.5. More precisely, for $v$ and $w$ words of length $p$ and 
respectively $q$ over the alphabet $[s]$, we get (by exactly the same 
argument as in the proof of Theorem 7.5) that:
\[
\lim_{k \to \infty} N_k^2 \Bigl( \ \E ( \, \tr 
( Y_v^{(k)} ) \cdot \tr ( Y_w^{(k)} ) \, ) -
\E ( \, \tr ( Y_v^{(k)} ) \, ) \cdot
\E ( \, \tr ( Y_w^{(k)} ) \, )  \ \Bigr) \ = 
\]
\begin{equation}
\card \Bigr\{ \tau \in \sanc \ \begin{array}{cl}
\vline & \tau \mbox { is a complete matching of } [p+q], \\
\vline & \tau \mbox { is $(p,q)$--connected, and } 
           (v \cup w) \circ \tau = v \cup w
\end{array} \Bigr\}  .
\end{equation}

The Equation (7.19) is supplementing the fact that 
\begin{equation}
\lim_{k \to \infty} \E \Bigl( \tr ( Y_v^{(k)} ) \cdot \tr ( Y_w^{(k)} )
\Bigr) \ = \ \lim_{k \to \infty} \E \Bigl( \tr ( Y_v^{(k)} ) \Bigr) \cdot
\lim_{k \to \infty} \E \Bigl( \tr ( Y_w^{(k)} ) \Bigr) ,
\end{equation}
where the limits on the right-hand side of (7.20) are described by using 
non-crossing permutations (or partitions) in the disc. The precise formula
for these limits (for instance for the word $w$ which has length $q$) is
\[
\lim_{k \to \infty} \E \Bigl( \tr ( Y_w^{(k)} ) \Bigr)  \ = \ 
\card \Bigr\{ \tau \in \GS_{disc-nc} ( \, [q] \, ) \ \begin{array}{cl}
\vline & \tau \mbox { is a complete matching of } [q] \\
\vline & \mbox { and } w \circ \tau = w
\end{array} \Bigr\}  ,
\]
and is well-known (in the language of free probability, this is the 
formulation in terms of non-crossing cumulants for the well-known fact that
$Y_1^{(k)}, \ldots , Y_s^{(k)}$ behave asymptotically like a semicircular
system of Voiculescu -- see \cite{V}).

$\ $

$\ $

\section{Non-crossing permutations in a multi-annulus}

{\bf 8.1 Notations.} In this section we fix a family of positive integers
\setcounter{equation}{0}
$p_1, \ldots , p_l$, where $l \geq 1$. For $1 \leq i \leq l$ we will denote 
\[
\gamma_i := (p_1+ \cdots + p_{i-1}+1, \ldots , p_1+ \cdots + p_{i-1}+p_i)
\in \GS ( \, [p_1 + \cdots + p_l] \, ).
\]
We denote $\gamma := \gamma_1 \gamma_2 \cdots \gamma_l$ (commuting product).

$\ $

{\bf 8.2 Remark.} Let $\tau$ be a permutation in 
$\GS ( \, [p_1 + \cdots + p_l] \, )$. Same as in Section 2.11,
we will denote as ``$\# ( \tau \vee \gamma )$'' the number of orbits into 
which $[p_1 + \cdots + p_l]$ is split under the joint action of $\tau$ and 
$\gamma$ (with $\gamma$ defined as above). Because of the pictorial 
interpretation presented in the next remark, we will refer to the situation
when $\# ( \tau \vee \gamma ) =1$ by saying that $\tau$ is 
{\em $(p_1, \ldots , p_l)$--connected.}

The inequality (2.10) of Section 2.11 gives us here that
\begin{equation}
\# ( \tau ) \ + \ \# ( \tau^{-1} \gamma ) \ \leq \
p_1 + \cdots + p_l -l + 2 \cdot \# ( \tau \vee \gamma ),
\end{equation}
for every $\tau \in \GS ( \, [p_1+ \cdots +p_l] \, )$.
This can also be re-written in terms of distances in the Cayley graph
of $\GS ( \, [ p_1 + \cdots + p_l ] \, )$, in the form:
\begin{equation}
d( \, id , \tau \, ) + d( \, \tau, \gamma \, ) \ \geq \ 
d( \, id , \gamma \, ) + 2 \Bigl( l- \# ( \tau \vee \gamma ) \Bigr).
\end{equation}
As such, it provides us with a lower bound on how close $\tau$ can 
be from lying on a geodesic between $id$ and $\gamma$, 
if $\# ( \tau \vee \gamma )$ is given.

$\ $

{\bf 8.3 Remark.} For the positive integers $p_1, \ldots , p_l$ fixed
in this section,  we consider the generalization of the pictorial 
framework used in Sections 3--6. We thus look at a system of $l$ circles: 
one ``external'' circle and $l-1$ ``internal'' ones, such that the closed 
discs enclosed by the internal circles are pairwise disjoint and are all 
contained in the open disc enclosed by the external circle. On the external
circle we mark $p_1$ points, labelled $1, \ldots , p_1$. Then we visit the 
$l-1$ internal circles (in some order), and mark on them $p_2, \ldots , p_l$
points (respectively), which we label with the integers from the intervals
$[ p_1 +1, \ldots , p_1 + p_2 ],$
$\ldots , [p_1 + \cdots +p_{l-1}+1, \ldots ,p_1+ \cdots +p_{l-1}+p_l],$
respectively. The points on the external circle are labelled in clockwise
order, while the points on each of the internal circles are labelled 
counter-clockwise. The multi-annulus comprised between the external circle 
and the $l-1$ internal circles will be referred to as ``{\em the 
$(p_1, \ldots , p_l)$--annulus}''. 

The discussion about planar pictures of permutations made in Section 3.1
generalizes verbatim to the case of the $(p_1, \ldots , p_l)$--annulus. The
only difference is that the planar maps which we draw have now $l$ (instead
of 2) red faces. (See Figure 11 for some concrete examples of admissible
and non-admissible drawings of a cycle in the $(3,3,3)$--annulus.) We will
use the name ``{\em $(p_1, \ldots , p_l)$--planar}'' for a permutation $\tau$ 
of $[ p_1 + \cdots + p_l ]$ with the property that admissible pictures of 
the cycles of $\tau$ can be drawn in the $(p_1, \ldots , p_l)$--annulus, 
such that the regions enclosed by these cycles are pairwise disjoint.

The algebraic counterpart for the concept of a $(p_1, \ldots , p_l)$--planar 
permutation should be the one of a ``{\em $(p_1, \ldots , p_l)$--annular
non-crossing permutation}''; this should be defined by identifying the 
crossing patterns that are to be avoided in the 
$(p_1, \ldots , p_l)$--annulus. It is not so clear what would be a a nice 
way of doing this (except for the case $l=2$ which is treated in the 
Section 3, and, of course, for $l=1$).

$\ $

{\bf 8.4 Problem.} Find an algebraic definition, involving crossing patterns,
which introduces the set $\smanc$ of $(p_1, \ldots , p_l)$--annular 
non-crossing permutations.

$\ $

The definition for the fact that $\tau \in \smanc$ should go via ``localized''
conditions, similar to those known in the cases when $l=1$ and $l=2$. If any 
extrapolation can be made based on these two cases, a localized condition for 
a given $l \geq 3$ should involve not more than $2l+2$ elements of 
$[p_1 + \cdots + p_l],$ belonging to not more than $l+1$ orbits of $\tau$.
Moreover, in the case when $\tau$ is a complete matching of 
$[ p_1+ \cdots + p_l ],$ the definition for ``$\tau \in \smanc$''
should reduce to the fact that any group of up to $l+1$ orbits ( = pairs) of
$\tau$ is $(p_1, \ldots , p_l)$--annular non-crossing.

The results proved in the Sections 5 and 6 for $l=2$ suggest some problems
one can pose in the $(p_1, \ldots , p_l)$--annular framework.

$\ $

{\bf 8.5 Problems.} Supposing that a suitable definition for $\smanc$ was 
found (in Problem 8.1), is it true that:

$1^o$ $\smanc$ is the smallest subset of $\GS ( \, [p_1 + \cdots + p_l] \, )$
which contains $\GS_{disc-nc} (p_1$
$+ \cdots + p_l )$ and is invariant under 
conjugation with $\gamma_1 , \ldots , \gamma_l$?

$2^o$ A permutation $\tau$ of $[p_1 + \cdots + p_l]$ belongs to $\smanc$ if 
and only if (8.1) holds with equality?

$\ $

{\bf 8.6 Remarks.} $1^o$ We are confident that the above question $2^o$ has 
an affirmative answer. This is because a suitable resolution of the 
Problem 8.4 will have to take into $\smanc$ precisely those permutations
of $[p_1 + \cdots + p_l]$ which are $(p_1, \ldots , p_l)$--planar (in the 
sense discussed in Remark 8.3); and it is quite plausible 
that a permutation $\tau$ of $[p_1+ \cdots +p_l]$ is 
$(p_1, \ldots , p_l)$--planar if and only if it satisfies (8.1) with 
equality. Indeed, if $\tau$ is $(p_1, \ldots , p_l)$--planar then let us look 
at the tri-coloured planar map giving its picture. This map has $l$ red faces
and $\# ( \tau )$ black faces. Under the believable assumption that 
Kreweras complementation does generalize to the 
$(p_1, \ldots , p_l)$--framework, we see moreover that the map has 
$\# ( \tau^{-1} \gamma )$ white faces. Hence the total number of faces of the 
map is $\# ( \tau ) + \# ( \tau^{-1} \gamma ) +l,$ and (as is immediately
checked) the Euler's characteristic formula amounts precisely to having 
equality in (8.1). For the converse, the idea is to create abstractly a 
compact connected surface by glueing together $l$ red faces, $\# ( \tau )$ 
black faces (which represent the cycles of $\tau$) and 
$\# ( \tau^{-1} \gamma )$ white faces (which represent 
the cycles of $\tau^{-1} \gamma$); the equality in (8.1) will imply that 
the surface so created has Euler characteristic equal to 2, hence 
is a sphere. But this sphere comes by construction with a
$(p_1, \ldots , p_l)$--planar picture of $\tau$ on it -- hence $\tau$ is
$(p_1, \ldots , p_l)$--planar.

$2^o$ In Problem 8.5.1, planarity considerations support the idea that 
we should have in any case the inclusion 
$\smanc \supset \GS_{disc-nc} (p_1 + \cdots + p_l).$ This is because if 
$\tau \in \GS_{disc-nc} (p_1 + \cdots + p_l),$ then a planar picture for 
$\tau$ in the disc can be ``folded'' to become a picture on an $l$--punctured 
sphere (a sphere with $l$ little circular punctures in it), and the latter
picture can in turn be transformed into a planar $(p_1, \ldots ,p_l)$--annular
picture for $\tau$. Thus it is quite likely that at least the inclusion 
``$\supset$'' in Problem 8.5.1 will have to be true (but it is less clear 
what to expect concerning the opposite inclusion ``$\subset$'').
 
$\ $

$\ $

\section{Asymptotic Gaussianity for traces of words made with 
Wishart matrices}

In this section we consider again the family of Wishart random matrices 
\setcounter{equation}{0}
$X_1, \ldots , X_s$ which appeared in Section 7 (with its version 
$X_1^{(k)}, \ldots , X_s^{(k)}$ used for discussing asymptotics -- cf.
Notations 7.3). We adopt all the notations introduced in the Section 7
in connection to traces of words $X_w$ formed with the matrices 
$X_1, \ldots , X_s$. We will now look at products of $l$ such traces, 
where we allow that $l \geq 3$. The computations are quite similar to those 
shown in Section 7 for the case $l=2$. Nevertheless, it is worth writing 
down the resulting formulas, due to an immediate consequence which they 
bear on a phenomenon of asymptotic Gaussianity (see Corollary 9.4 below). 
It is interesting to observe that in order to obtain this corollary on 
asymptotic Gaussianity, one does not need to make any kind of analysis 
of multi-annular non-crossing permutations for $l \geq 3$; indeed, for 
$l \geq 3$ the general inequality (8.1) is all that is needed (the case 
when (8.1) holds with equality does not appear in the discussion).

Let us now elaborate. The first thing to do is generalize the formulas 
presented in Lemma 7.6. We leave it as an exercise to the reader to check 
that the proof of Lemma 7.6 (shown in the case $l=2$) extends mutatis 
mutandis to give us the following:

$\ $

{\bf 9.1 Lemma.} Let $w_1, \ldots , w_l$ be words of lengths
$p_1, \ldots , p_l$ (respectively)  over the alphabet $[s]$, where 
$l \geq 1$ and $p_1, \ldots , p_l \geq 1$. Let $w$ be the word 
of length $p_1 + \cdots + p_l$ which is obtained by juxtaposing 
$w_1, \ldots , w_l$ (in this order), and let the permutation 
$\gamma \in \GS ( \, [p_1 + \cdots + p_l] \, )$ be as in the 
Notations 8.1. Then 
\begin{equation}
\E \Bigl( \tr (X_{w_1}) \cdots \tr (X_{w_l}) \Bigr) \ = \ 
\frac{1}{N^{p_1 + \cdots + p_l +l}} \sum_{ \begin{array}{c}
{\scriptstyle \tau \in \GS ( [p_1 + \cdots + p_l] ) } \\
{\scriptstyle such \ that \ w \circ \tau = w}
\end{array} } \ M^{ \# ( \tau ) } N^{ \# ( \tau^{-1} \gamma ) }.
\end{equation}

$\ $

The next point is to adjust the Equation (9.1) so that on its right-hand 
side we are only left with a summation over the permutations 
$\tau \in \GS ( \, [p_1 + \cdots + p_l] \, )$ which are 
$(p_1, \ldots , p_l)$--connected. This has to be a generalization for 
the Equation (7.16) which appeared in the proof of Theorem 7.5, in the 
case $l=2$. In that case the only thing which needed to be done was to 
take the variance $\E ( \, \tr (X_{w_1}) \cdot \tr (X_{w_2}) \, ) -
\E ( \, \tr (X_{w_1}) \, ) \cdot \E ( \, \tr (X_{w_2}) \, )$. For general 
$l$ as we have in (9.1), one needs to form some more complicated expressions, 
called the ``cumulants'' of the random variables 
$\tr ( X_{w_1} ), \ldots , \tr (X_{w_l})$.

$\ $

{\bf 9.2 Review} {\em of cumulants.} For every $n \geq 1$, the 
{\em cumulant of order $n$} is a certain multilinear functional 
$\Cl_n : \Bigl( \ L^{\infty}_{-} ( \Omega , \F, P) \ \Bigr)^n \to \C$,
where $L^{\infty}_{-} ( \Omega , \F, P)$ is as introduced 
in Section 7.1. For instance for $n \leq 3$ the formulas defining 
$\Cl_n$ are as follows:
\begin{equation}
\left\{  \begin{array}{ccl}
\Cl_1 ( g_1)          & =  &   \E (g_1)                                  \\
\Cl_2( g_1, g_2)      & =  &   \E (g_1 g_2) - \E ( g_1 ) \cdot \E (g_2)  \\
\Cl_3( g_1, g_2, g_3) & =  &
                        \E (g_1 g_2 g_3) - \E ( g_1 ) \cdot \E (g_2 g_3) 
                                   - \E ( g_2 ) \cdot \E( g_1 g_3 )      \\
                      &     &   - \E (g_1 g_2 ) \cdot \E (g_3) 
                             + 2 \E ( g_1 ) \cdot \E ( g_2 ) \cdot \E (g_3 )
\end{array} \right.
\end{equation}
(for $g_1, g_2, g_3 \in L^{\infty}_{-} ( \Omega , \F, P)$).

In order to give explicitly the formula defining $\Cl_n$ for an arbitrary 
$n \geq 1$, one uses the Moebius function for set-partitions, as 
described for instance in \cite{R}. For every $n \geq 1$, let us denote
by $\PP_n$ the poset of all partitions of $[n]$, where the partial order 
on $\PP_n$ is defined by setting $\pi \leq \rho$ if and only if every block 
of $\pi$ is contained in a block of $\rho$. The Moebius function for 
this poset is a function 
$$\mu_n : \{ ( \pi , \rho ) \ | \ \pi , \rho \in \PP_n, \ 
\pi \leq \rho \} \to \Z ,$$
uniquely determined by the fact that it satisfies
\begin{equation}
\left\{  \begin{array}{ccll}
\mu_n ( \pi , \pi ) & =  & 1, & \forall \ \pi \in \PP_n    \\
\begin{array}[t]{c} \sum \\ \pi \in [ \theta , \rho ] \end{array} \
\mu_n ( \pi , \rho ) & = & 0 &
\forall \ \theta , \rho \in \PP_n \mbox{ such that }  
\theta \leq \rho , \ \theta \neq \rho         
\end{array}  \right.
\end{equation}
(where $[\theta , \rho ] := \{ \pi \in \PP_n \ | \ \theta \leq 
\pi \leq \rho \}$). The meaning of the Equations (9.3) is that they make 
$\mu_n$ become the inverse of the function identically equal to 1, under 
a certain convolution operation (see \cite{R}, Section 3). 

The explicit formula defining the cumulant functional $\Cl_n$ is then:
\begin{equation}
\Cl_n ( g_1, \ldots , g_n ) = \sum_{\pi \in \PP_n} 
\E_{\pi} (g_1, \ldots , g_n) \mu_n ( \pi , 1_n)
\end{equation}
for $n \geq 1$ and $g_1, \ldots , g_n \in L^{\infty}_{-} ( \Omega , \F, P)$,
where $1_n \in \PP_n$ is the partition of $[n]$ with only one block, and where 
for $\pi = \{ B_1 , \ldots , B_m \} \in \PP_n$ we denote
\begin{equation}
\E_{\pi} ( g_1 , \ldots , g_n ) \ = \ 
\E ( \  \prod_{i \in B_1} g_i \ ) \cdots
\E ( \  \prod_{i \in B_m} g_i \ ) .
\end{equation}
(The reader not used to these notations could practice them by verifying that 
for $n = 1,2,3$ the Equation (9.4) reduces indeed to what had been announced
in (9.2). Concerning concrete computations for the numbers 
$\mu_n ( \pi , \rho )$: it is fairly easy to find a general 
explicit formula for them -- see e.g. the corollary at the end of Section 7
in \cite{R}.)

The importance of the cumulant functionals comes from the following fact:
the generating series for the cumulants of a family of random variables is
essentially the log of the characteristic function of the family (see
e.g. \cite{Sh}, Section 12 in Chapter II). In particular, we have the 
following Gaussianity criterion in terms of cumulants: a family 
$\{ g_{\lambda} \ | \lambda \in \Lambda \}$ of random variables in 
$L^{\infty}_{-} ( \Omega , \F, P)$ is Gaussian if and only if the 
cumulants $\CC_n ( g_{\lambda_1}, \ldots , g_{\lambda_n} )$ vanish for
all choices of $n \geq 3$ and of $\lambda_1 , \ldots , \lambda_n \in 
\Lambda$.

$\ $

{\bf 9.3 Proposition.} In the framework of Lemma 9.1, we have that 
\begin{equation}
\CC_l \Bigl( \tr (X_{w_1}) , \ldots , \tr (X_{w_l}) \Bigr) \ = \ 
\frac{1}{N^{p_1 + \cdots + p_l +l}} \sum_{ \begin{array}{c}
{\scriptstyle \tau \in \GS ( [p_1 + \cdots + p_l] ) } \\
{\scriptstyle such \ that \ w \circ \tau = w}  \\
{\scriptstyle and \ \tau \ is \ (p_1, \ldots , p_l)-connected}
\end{array} } \ M^{ \# ( \tau ) } N^{ \# ( \tau^{-1} \gamma ) }.
\end{equation}

$\ $

{\bf Proof.} Let us denote the orbits of $\gamma$ by $J_1, \ldots , J_l$,
where $J_1 = \{ 1, \ldots , p_1 \}$,
$J_2 = \{ p_1 + 1, \ldots , p_1 + p_2 \} , \ldots , J_l$ =
$\{ p_1 + \cdots + p_{l-1} + 1, \ldots , p_1 + \cdots + p_l \}$.
A permutation $\tau \in \GS ( [ p_1 + \cdots + p_l ] )$ induces a 
partition $\alpha ( \tau )$ of $[l]$, via the requirement that 
$i,j \in [l]$ are in the same block of $\alpha ( \tau )$ if and only if 
$J_i$ and $J_j$ are contained in the same orbit of the joint action of 
$\tau$ and $\gamma$ on $[n]$. (The definition of $\alpha ( \tau )$ makes 
sense because every orbit of the joint action of $\tau$ and $\gamma$ is a 
union of $J_i$'s.) Note that in particular a permutation 
$\tau \in \GS ( [ p_1 + \cdots + p_l ] )$ is $(p_1, \ldots , p_l)$-connected
if and only if $\alpha ( \tau )$ is equal to $1_l$, the partition of $[l]$
which has only one block.

The induced partitions ``$\alpha ( \tau )$'' from the preceding paragraph 
play a role when one takes the argument which led to Equation (7.15) in 
the proof of Theorem 7.5, and adjusts it to work for a more complicated 
product of traces of words. 
[Concrete example: if we suppose just for a moment that $l=3$, and work out
the analogue of Equation (7.15) for
$\E ( \, \tr (X_{w_1} ) \cdot \tr ( X_{w_3}) \, ) \cdot \E ( \tr (X_{w_2} ))$,
then it is immediate that we will get a sum indexed by permutations 
$\tau \in \GS ( [p_1+p_2+p_3 ] )$ which satisfy $w \circ \tau = w$ and the 
extra condition $\alpha ( \tau ) \leq \{ \, \{ 1,3 \} , \{ 2 \} \, \}$ in 
$\PP_3$.]
The reader should have no difficulty to verify that the resulting formula 
can be stated in general 
as follows: for every partition $\pi \in \PP_l$, we have
\[
\E_{\pi} \Bigl( \, \tr ( X_{w_1} ), \ldots , \tr ( X_{w_l} ) \, \Bigr) 
\ = \ \frac{1}{N^{p_1 + \cdots + p_l +l}} \sum_{ \begin{array}{c}
{\scriptstyle \tau \in \GS ( [p_1 + \cdots + p_l] ) } \\
{\scriptstyle such \ that \ w \circ \tau = w}  \\
{\scriptstyle and \ \alpha( \tau ) \leq \pi}
\end{array} } \ M^{ \# ( \tau ) } N^{ \# ( \tau^{-1} \gamma ) }
\]
(where the expression $\E_{\pi} ( \ \cdots \ )$ is as introduced in 
Equation (9.5) above). It is convenient to write this in a more 
compressed form as 
\begin{equation}
\E_{\pi} \Bigl( \, \tr ( X_{w_1} ), \ldots , \tr ( X_{w_l} ) \, \Bigr) 
\ = \ \sum_{\theta \leq \pi \ in \ \PP_l} \ Q_{\theta},
\end{equation}
where for every $\theta \in \PP_l$ we set
\begin{equation}
Q_{\theta} \ := \ 
\frac{1}{N^{p_1 + \cdots + p_l +l}} \sum_{ \begin{array}{c}
{\scriptstyle \tau \in \GS ( [p_1 + \cdots + p_l] ) } \\
{\scriptstyle such \ that \ w \circ \tau = w}  \\
{\scriptstyle and \ \alpha ( \tau ) = \theta}
\end{array} } \ M^{ \# ( \tau ) } N^{ \# ( \tau^{-1} \gamma ) } .
\end{equation}
But then, by substituting the Equation (9.7) in (9.4) we get that
\begin{equation}
\CC_l \Bigl( \ \tr ( X_{w_1} ), \ldots , \tr ( X_{w_l} ) \ \Bigr)
\ = \ \sum_{\pi \in \PP_l} \ \Bigl( \ 
\sum_{\theta \leq \pi \ in \ \PP_l} \ Q_{\theta} \ \Bigr) \
\mu_l ( \pi , 1_l).
\end{equation}
When we interchange the order of summation, the right-hand side of
(9.9) becomes 
\[
\sum_{\theta \in \PP_l} \ Q_{\theta} \cdot \Bigl( \
\sum_{\pi \geq \theta \ in \PP_l} \ \mu_l ( \pi , 1_l ) \ \Bigr)
\]
which is in turn equal to just $Q_{1_l}$, by virtue of the Equations (9.3)
satisfied by the Moebius function. So in conclusion we obtained that
\[
\CC_l \Bigl( \ \tr ( X_{w_1} ), \ldots , \tr ( X_{w_l} ) \ \Bigr)
\ = \ Q_{1_l},
\]
which is precisely (9.6) (by the definition of $Q_{1_l}$ in (9.8) and 
the fact that the condition ``$\alpha ( \tau ) = 1_l$'' is equivalent
to the $(p_1, \ldots , p_l)$-connectedness of $\tau$).   {\bf QED}

$\ $

We now arrive at the corollary on asymptotic Gaussianity which was 
announced at the beginning of the section. In this corollary we will
revert to the notations with an extra index $k$, and where 
$M_k / N_k \to c \in (0, \infty )$ as $k \to \infty$. Moreover, in the
Corollary 9.4 we will actually assume that the sequences 
$(M_k)_{k=1}^{\infty}$ and $(N_k)_{k=1}^{\infty}$ are picked in such a 
way that the stronger limit condition 
\begin{equation}
\lim_{k \to \infty} M_k - cN_k = c'
\end{equation}
is holding (for some $c' \in \R$).

$\ $

{\bf 9.4 Corollary.} Consider the hypotheses described above, and for every 
word $w$ of length $n$ over the alphabet $[s]$, consider the quantities
\begin{equation}
\EE_w \ := \sum_{  \begin{array}{c}
{\scriptstyle \tau \in \sdnc }  \\
{\scriptstyle such \ that \ w \circ \tau = w}
\end{array} } \ c^{\# ( \tau )}  \ \mbox{ and } \ 
\EE_w ' \ := \sum_{  \begin{array}{c}
{\scriptstyle \tau \in \sdnc }  \\
{\scriptstyle such \ that \ w \circ \tau = w}
\end{array} } \ \# ( \tau ) \cdot c^{\# ( \tau ) -1 } \cdot c' .
\end{equation}
Then the family of random variables 
\begin{equation}
\left\{ N_k  \tr (X_w^{(k)} ) - N_k \EE_w - \EE_w '\ | \ 
\mbox{$w$ word over the alphabet $[s]$} \right\}
\end{equation}
is asymptotically a centered Gaussian family for $k \to \infty$.

$\ $

{\bf Proof.} We start by explaining why, for a given word $w$, the random 
variable $N_k   \tr (X_w^{(k)} ) - N_k \EE_w - \EE_w '$ is asymptotically
centered for $k \to \infty$. Let us write $\E ( \tr ( X_w^{(k)} )$ as in the 
Lemma 7.6.1,
\[
\E ( \tr ( X_w^{(k)} ) \ = \ 
\frac{1}{N_k^{n+1}} \sum_{ \begin{array}{c}
{\scriptstyle \tau \in \GS ( [n]) \ such}  \\
{\scriptstyle that \ w \circ \tau = w}
\end{array} } \ M_k^{\# ( \tau )} N_k^{\# ( \tau^{-1} \gamma_o)}
\]
\begin{equation}
= \ \sum_{ \begin{array}{c}
{\scriptstyle \tau \in \GS ( [n]) \ such}  \\
{\scriptstyle that \ w \circ \tau = w}
\end{array} } \ \Bigl( \frac{M_k}{N_k} \Bigr)^{\# ( \tau )} 
\cdot N_k^{( \# ( \tau ) + \# ( \tau^{-1} \gamma_o )) - (n+1)}
\end{equation}
(where $n$ is the length of $w$ and $\gamma_o = (1,2, \ldots , n) \in 
\GS ( [n] )$). We know that the exponent 
$( \# ( \tau ) + \# ( \tau^{-1} \gamma_o )) - (n+1)$ appearing in (9.13)
is always $\leq 0$, with equality if and only $\tau \in \GS_{disc-nc} (n)$.
The additional point to be noticed here is that (due to an easy argument 
on signatures, which is left to the reader), the difference 
$( \# ( \tau ) + \# ( \tau^{-1} \gamma_o )) - (n+1)$ is even for every
$\tau \in \GS ( [n] )$, and it is hence $\leq -2$ whenever $\tau$ is not
in $\GS_{disc-nc} ( [n] )$. Consequently, the terms of the sum in (9.13) 
which are indexed by permutations $\tau \in \GS ( [n] ) \setminus \sdnc$ 
will converge to 0 even after being multiplied by $N_k$. The rest of the 
sum in (9.13) is
\[
\EE_w^{(k)} \ := \ \sum_{ \begin{array}{c}
{\scriptstyle \tau \in \sdnc \ such}  \\
{\scriptstyle that \ w \circ \tau = w}
\end{array} } \ \Bigl( \frac{M_k}{N_k} \Bigr)^{\# ( \tau )} ,
\]
and we are left with showing that 
\[
\lim_{k \to \infty} N_k ( \EE_w^{(k)} - \EE_w ) - \EE_w ' \ = \ 0.
\]
And indeed, from the definition of $\EE_w$ and $\EE_w^{(k)}$ we have that
\[
N_k ( \EE_w^{(k)} - \EE_w ) \ = \ \sum_{ \begin{array}{c}
{\scriptstyle \tau \in \sdnc \ such}  \\
{\scriptstyle that \ w \circ \tau = w}
\end{array} } \ N_k \cdot \Bigl( \ \Bigl( \frac{M_k}{N_k} \Bigr)^{\# ( \tau )}
- c^{\# ( \tau )} \  \Bigr)
\]
\[
= \ \sum_{ \begin{array}{c}
{\scriptstyle \tau \in \sdnc \ such}  \\
{\scriptstyle that \ w \circ \tau = w}
\end{array} } \ (M_k - c N_k)  \cdot \Bigl( \  \sum_{j=0}^{\# ( \tau ) -1}
\Bigl( \frac{M_k}{N_k} \Bigr)^j \cdot c^{\# ( \tau ) -1-j} \Bigr),
\]
and the latter quantity clearly converges to $\EE_w '$ for $k \to \infty$.

We now take on the asymptotic Gaussianity of the family in (9.12). We will
prove it by verifying the asymptotic vanishing of all the cumulants of order 
$l \geq 3$ made with random variables from the family (cf. review in 
Section 9.2). For any $l \geq 3$ and any words $w_1, \ldots , w_l$ over the 
alphabet $[s]$ we have:
\[
\CC_l \left( \ 
N_k \tr (X_{w_1}^{(k)}) - N_k \EE_{w_1} - \EE_{w_1} ' , \ldots ,
N_k \tr (X_{w_l}^{(k)}) - N_k \EE_{w_l} - \EE_{w_l} '   \ \right)
\]
\[
= \ \CC_l \Bigl( \ 
N_k \tr (X_{w_1}^{(k)}) , \ldots , N_k \tr (X_{w_l}^{(k)}) \ \Bigr) 
\]
(because a cumulant of order $\geq 2$ does not change when some constants 
are added to its arguments)
\[
= \ N_k^l \cdot \CC_l \Bigl( \ \tr (X_{w_1}^{(k)}) , \ldots ,
\tr (X_{w_l}^{(k)}) \ \Bigr) 
\ \ \mbox{ (by multi-linearity)}
\]
\begin{equation}
= \ \frac{N_k^l}{N^{p_1 + \cdots + p_l +l}} \sum_{ \begin{array}{c}
{\scriptstyle \tau \in \GS ( [p_1 + \cdots + p_l] ) } \\
{\scriptstyle such \ that \ w \circ \tau = w}  \\
{\scriptstyle and \ \tau \ is \ (p_1, \ldots , p_l)-connected}
\end{array} } \ M_k^{ \# ( \tau ) } N_k^{ \# ( \tau^{-1} \gamma ) } ,
\end{equation}
where in (9.14) we denoted the lengths of $w_1, \ldots , w_l$ by 
$p_1, \ldots , p_l$, and went into the framework and notations of 
Proposition 9.3. The quantity in (9.14) can be further written as 
\begin{equation}
\sum_{ \begin{array}{c}
{\scriptstyle \tau \in \GS ( [p_1 + \cdots + p_l] ) } \\
{\scriptstyle such \ that \ w \circ \tau = w}  \\
{\scriptstyle and \ \tau \ is \ (p_1, \ldots , p_l)-connected}
\end{array} } \ \Bigl( \frac{M_k}{N_k} \Bigr)^{ \# ( \tau ) } \cdot
N_k^{ ( \# ( \tau ) + \# ( \tau^{-1} \gamma )) - (p_1 + \cdots +p_l) }.
\end{equation}
But now, for any $(p_1, \ldots , p_l)$--connected permutation $\tau$, 
the basic inequality (8.1) (used with $\# ( \tau \vee \gamma ) = 1$)
gives us that
\[
\# ( \tau ) + \# ( \tau^{-1} \gamma ) \leq ( p_1 + \cdots + p_l ) -l+2.
\]
This implies that in all the terms in the sum (9.15) the exponent of 
$N_k$ is at most $2-l$ (which is $\leq -1$), and the desired convergence 
to 0 immediately follows.  {\bf QED}

$\ $

{\bf 9.5 Remark.} The Theorem 7.5 in Section 7 can now be interpreted 
as giving us a formula for the asymptotic covariance between two 
random variables (indexed by words $v,w$ of length $p,q$) in the 
family (9.12). We mention that the combinatorics of the sets 
$\sanc$ can be further used to understand how the asymptotically 
Gaussian family (9.12) should be transformed in order to also become 
asymptotically independent; this will be presented in the paper
\cite{KMS}. In the case when we have only one Wishart matrix 
(i.e. when $s=1$) both the statement of the Corollary 9.4 and the 
further discussion on how to obtain asymptotic independence were
derived by Cabanal-Duvillard \cite{C} by different methods (based 
on stochastic integrals rather than on combinatorics).

$\ $

$\ $

\setcounter{section}{8}\section*{Appendix}


\noindent\hbox to \hsize{\quad
\vtop{\hsize150pt\noindent
\includegraphics{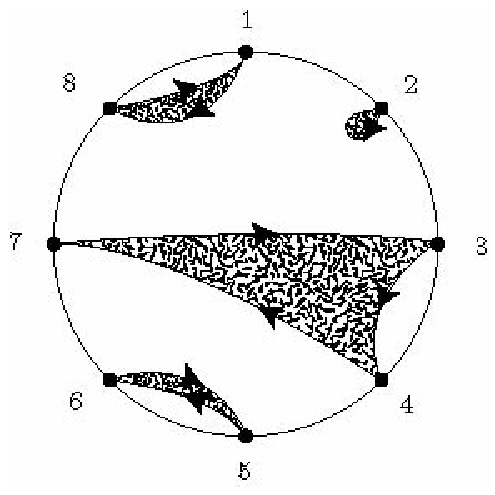}

\medskip

\capt{%
{\bf Figure 1.} $\tau = (1,8) (2) (3,4,7)\allowbreak
(5,6) \allowbreak \in {\cal S}_{\hbox{{\scriptsize
\it disc-nc}}} (8)$}}%
\hfill
\vtop{\hsize150pt\noindent
\includegraphics{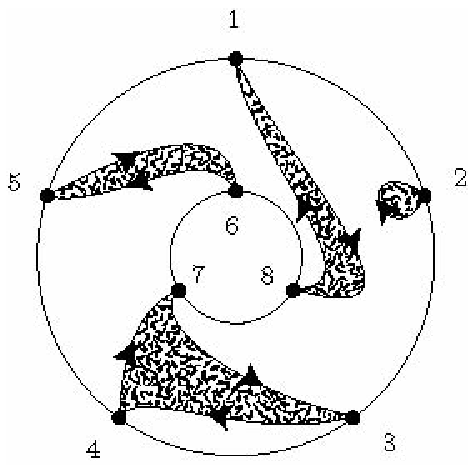}

\medskip

\capt{%
{\bf Figure 2.} $\tau = (1,8) (2) (3,4,7)\allowbreak
(5,6) \in {\cal S}_{\hbox{{\scriptsize\it ann-nc}}} (5,3)$}}\quad}


\bigskip

\noindent\hbox to \hsize{\quad
\vtop{\hsize150pt\noindent
\includegraphics{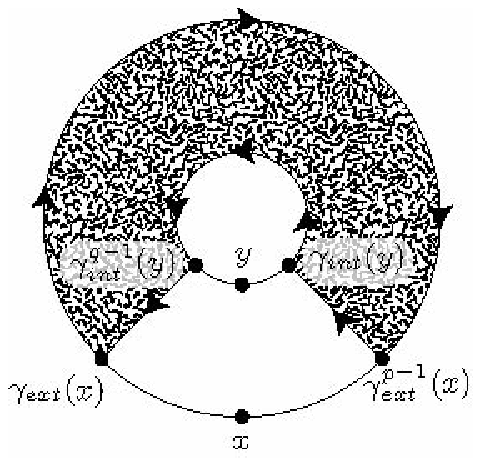}

\medskip

\capt{%
{\bf Figure 3.} The picture of $\lambda_{x,y}$}}%
\hfill
\vtop{\hsize150pt\noindent
\includegraphics{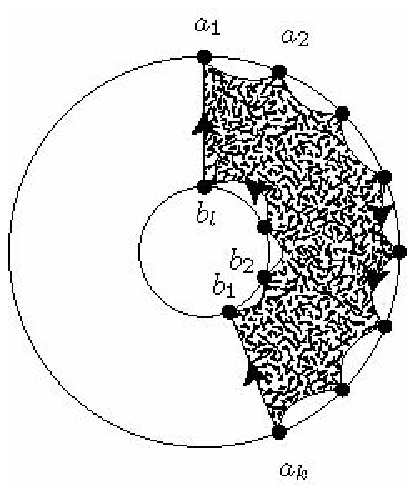}

\medskip

\capt{%
{\bf Figure 4.} A $(p,q)$-connecting cycle which is standard in
the annular sense}}\quad}


\bigskip

\noindent\hbox to \hsize{\qquad
\vtop{\hsize 400pt\noindent
\includegraphics{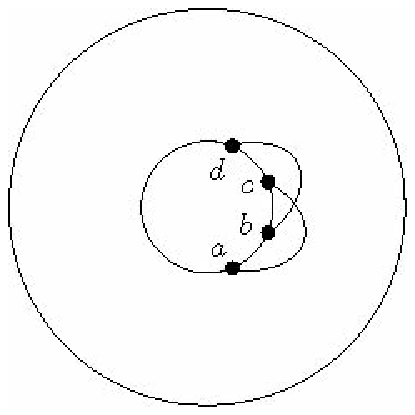}\hfill
\includegraphics{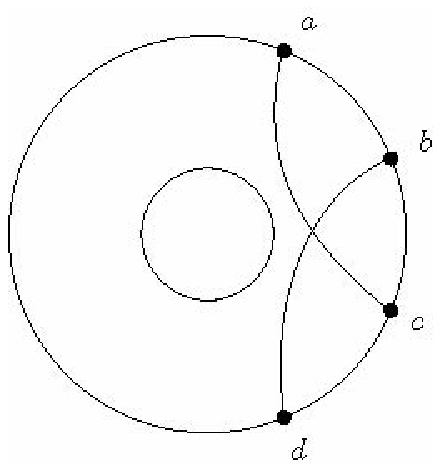}\qquad\hbox{}

\medskip

\hbox to \hsize{\hfill\small
{\bf Figure 5} Illustrations of the pattern (AC-1). \hfill}}}


\bigskip

\noindent\hbox to \hsize{\qquad
\vtop{\hsize 400pt\noindent
\includegraphics{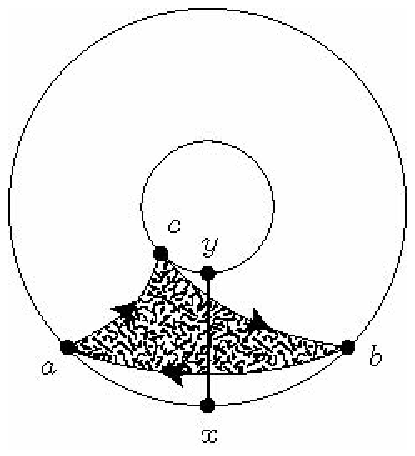}\hfill
\includegraphics{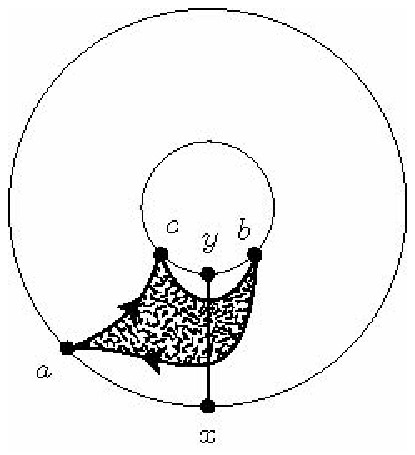}\qquad\hbox{}

\medskip

\hbox to \hsize{\hfill\small
{\bf Figure 6} Illustrations of the pattern (AC-2). \hfill}}}


\bigskip

\noindent\hbox to \hsize{\qquad
\vtop{\hsize 400pt\noindent
\includegraphics{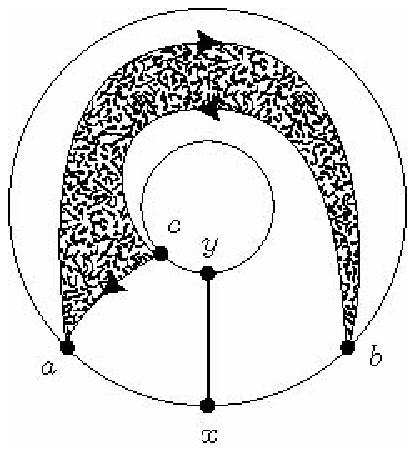}\hfill
\includegraphics{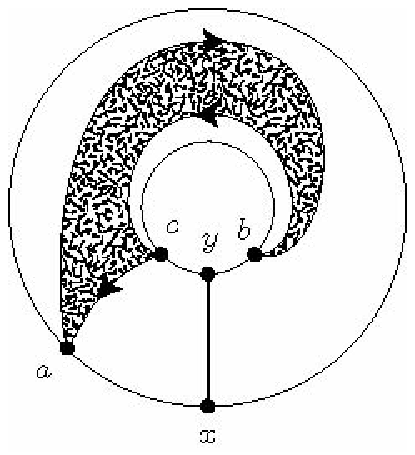}\qquad\hbox{}

\medskip

\hbox to \hsize{\hfill\doublecapt{%
{\bf Figure 6'.} The orientation of $\tau$ on $\{a, b, c\}$ is
essential in (AC-2). Unlike in figure 6, $\tau$ and $\lambda_{x,y}$
induce the same cycle on $\{ a, b, c\}$.} \hfill}}}


\bigskip

\noindent\hbox to \hsize{\qquad
\vtop{\hsize 450pt\noindent
\includegraphics{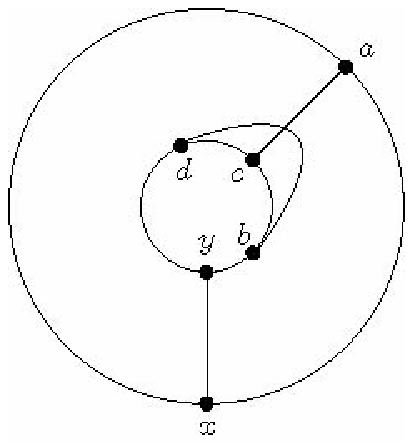}\quad\
\includegraphics{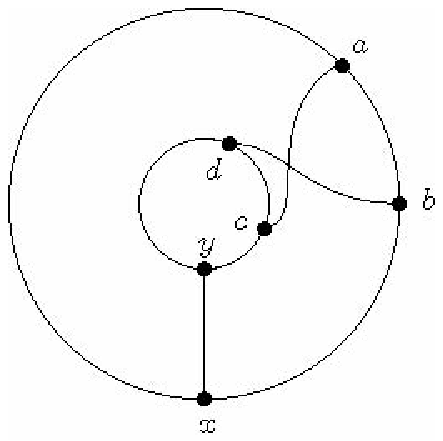}\quad
\includegraphics{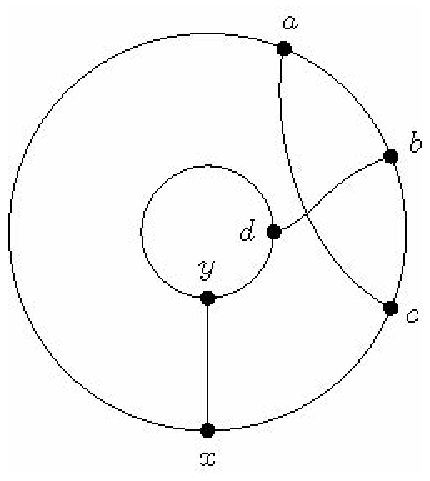}

\medskip

\hbox to \hsize{\hfill\doublecapt{%
{\bf Figure 7.} Illustrations of the crossing pattern (AC-3).}\hfill}}}


\bigskip

\noindent\hbox to \hsize{\quad
\vtop{\hsize150pt\noindent
\includegraphics{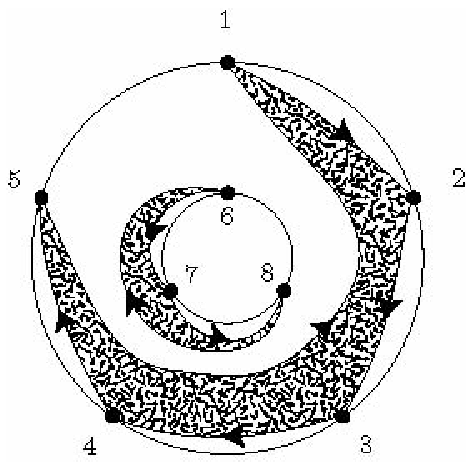}

\medskip

\capt{%
{\bf Figure 8.} $p= 5$, $q = 3$. Planar $(p,q)$-annular drawing for $\tau =
\gamma = (1,2,3,4,5) (6,7,8)$.}}\hfill
\vtop{\hsize150pt\noindent
\includegraphics{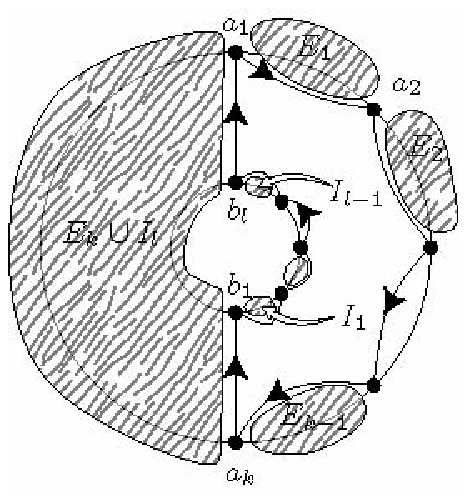}

\medskip

\capt{%
{\bf Figure 9.} Illustration for the discussion in Remark 3.9.}} \quad}


\bigskip

\noindent\hbox to \hsize{\quad
\vtop{\hsize150pt\noindent
\includegraphics{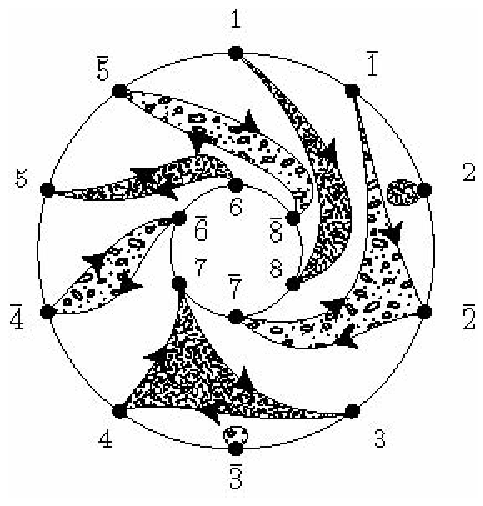}

\medskip

\capt{%
{\bf Figure 10.} The Kreweras complement of $\tau = (1,8) (2)
(3,4,7) (5,6) \in {\cal S}_{\hbox{{\scriptsize\it ann-nc}}} (5,3)$ is
$\tau^{-1}\gamma = (1, 2, 7) (3) (4, 6) (5, 8)$.}}%
\hfill
\vtop{\hsize150pt\noindent
\includegraphics{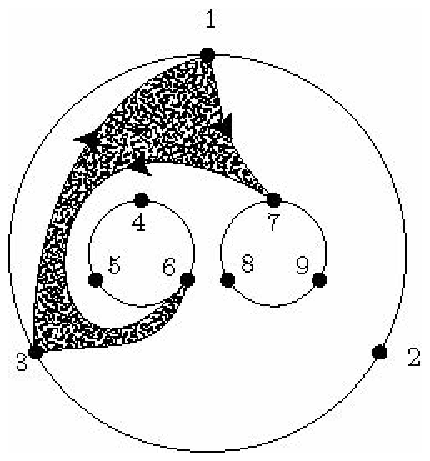}

\medskip

\capt{{\bf Figure 11a.}
An admissible drawing of the cycle
$(1,7,6,3)$ in the $(3,3,3)$ annulus. 
}}\quad}


\medskip

\noindent\hbox to \hsize{\qquad
\vtop{\hsize 400pt\noindent
\includegraphics{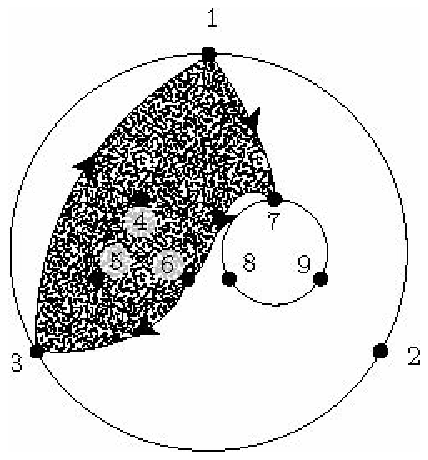}\hfill
\includegraphics{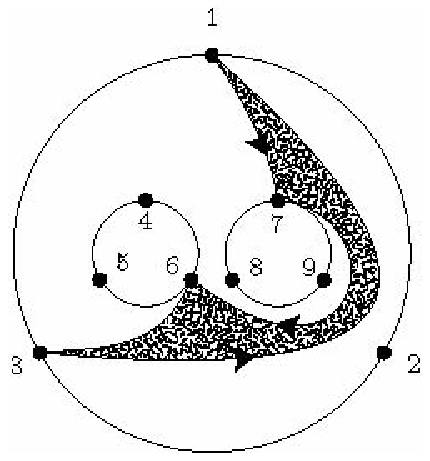}\qquad\hbox{}

      \medskip

\hbox to \hsize{\hfill\doublecapt{%
{\bf Figure 11b.} Inadmissible drawings of the cycle
$(1,7,6,3)$ in the $(3,3,3)$ annulus.}\hfill}}}

\setbox1=\vtop{\footnotesize\hsize 170pt\noindent
Alexandru Nica\\
Department of Pure Mathematics\\ 
University of Waterloo \\
Waterloo, Ontario N2L\,3G1, Canada            

\smallskip\noindent{\tt anica@math.uwaterloo.ca}}

\setbox2=\vtop{\footnotesize\hsize 170pt\noindent
James A. Mingo\\
Department of Mathematics and Statistics \\
Queen's University \\
Kingston, Ontario, K7L\,3N6, Canada 

\smallskip
\noindent
{\tt mingoj@mast.queensu.ca}}

\bigskip

\hbox to \hsize{\noindent\box2\hfill\box1}

\end{document}